\documentclass[12pt,a4paper,reqno]{amsart}
\usepackage{amsmath}
\usepackage{amsfonts}
\usepackage{amssymb}
\numberwithin{equation}{section}

\theoremstyle{plain}

\newtheorem{theorem}[subsection]{Theorem}
\newtheorem{proposition}[subsection]{Proposition}
\newtheorem{lemma}[subsection]{Lemma}
\newtheorem{corollary}[subsection]{Corollary}

\newtheorem{remark}[subsection]{Remark}
\newtheorem{example}[subsection]{Example}

\theoremstyle{definition}

\newtheorem{definition}[subsection]{Definition}

\renewcommand{\leq}{\leqslant}
\renewcommand{\geq}{\geqslant}

\newsavebox{\proofbox}
\savebox{\proofbox}{\begin{picture}(7,7)%
  \put(0,0){\framebox(7,7){}}\end{picture}}

\def\tr{\hbox{\rm tr}}

\def\B{{\mathcal B}}

\def\Z{{\mathbb Z}}
\def\E{{\mathbb E}}
\def\C{{\mathbb C}}
\def\R{{\mathbb R}}

\def\Q{{\mathbb Q}}
\def\P{{\mathbb P}}

\def\eps{{\varepsilon}}

\def\emph#1{{\it #1}}
\def\textbf#1{{\bf #1}}

\begin{document}

\title[Arithmetic progressions and primes]{Arithmetic progressions and the primes - El Escorial lectures}

\author{Terence Tao}
\address{Department of Mathematics\\University of California at Los Angeles\\ Los Angeles CA 90095}

\email{tao@math.ucla.edu}

\thanks{The author is
supported by a grant from the Packard Foundation.}


\begin{abstract}
We describe some of the machinery behind recent progress in establishing infinitely many arithmetic progressions of length $k$ in various sets of integers, in particular in arbitrary dense subsets of the integers, and in the primes.  
\end{abstract}

\maketitle

\section{Introduction}

A celebrated theorem of Roth \cite{roth} in 1953 asserts:

\begin{theorem}[Roth's theorem, first version]\cite{roth}
Let $A \subset \Z^+$ be a subset of integers with positive upper density, thus $\limsup_{N \to \infty} \frac{1}{N} |A \cap [1,N]| > 0$.
Then $A$ contains infinitely many arithmetic progressions $n, n+r, n+2r$ of length three.
\end{theorem}

Here we of course restrict the spacing of the progression $r$ to be non-negative.  This theorem was originally proven by Roth by Fourier analytic methods and a stopping time argment, and we shall reprove it below (in fact, we shall give two proofs).  This theorem was then generalized substantially by Szemer\'edi in 1975:

\begin{theorem}[Szemer\'edi's theorem, first version]\cite{szemeredi-4}, \cite{szemeredi}
Let $A \subset \Z^+$ be a subset of integers with positive upper density, thus $\limsup_{N \to \infty} \frac{1}{N} |A \cap [1,N]| > 0$,
and let $k \geq 3$.
Then $A$ contains infinitely many arithmetic progressions $n, n+r, \ldots, n+(k-1)r$ of length $k$.
\end{theorem}

Thus Roth's theorem is the $k=3$ version of Szemer\'edi's theorem.  (The cases $k < 3$ are trivial).

Szemer\'edi's original proof was combinatorial (relying in particular on graph theory) and very complicated.  A substantially shorter 
proof - but one involving the full machinery of
measure theory and ergodic theory, as well as the axiom of choice - was obtained by Furstenberg \cite{furst}, \cite{FKO} in 1977.
Since then, there have been two other types of proofs; a proof of Gowers \cite{gowers-4}, \cite{gowers} in 2001 which combines 
``higher order'' Fourier analytic methods with techniques from additive combinatorics; and also arguments of Gowers \cite{gowers-reg} and
Rodl-Skokan \cite{rodl}, \cite{rodl2} using the machinery of hypergraphs.  While we will not discuss all these separate proofs in detail here, we will need to discuss certain ideas
from each of these arguments as they will eventually be used in the proof of Theorem \ref{gt-thm} below.

The above theorems do not apply directly to the set of prime numbers, as they have density zero.  Nevertheless, in 1939 van der Corput \cite{van-der-corput} proved, by using Fourier analytic methods (the Hardy-Littlewood circle method) which were somewhat similar to the methods used by Roth, the following result:

\begin{theorem}[Van der Corput's theorem]\cite{van-der-corput}
Let $P \subset \Z^+$ be the set of primes. Then $P$ contains infinitely many arithmetic progressions $n, n+r, n+2r$ of length three.
\end{theorem}

However, just as Roth's Fourier-analytic methods proved very difficult to extend beyond the $k=3$ case, so too did van der Corput's arguments.
The proof relied on very delicate information concerning the Fourier coefficients of the primes (or more precisely of the von Mangoldt function
$\Lambda(n)$, which is essentially supported on the primes).  This additional information allows one to not only show that there are infinitely
many progressions of primes of length three, but also to obtain an asymptotic count as to \emph{how many} such progressions there are; we shall
return to this point later. 

Roth's theorem and van der Corput's theorem were combined by Green \cite{green} in 2003 to obtain

\begin{theorem}[Green's theorem]\label{gthm}\cite{green}
Let $A \subset P$ be a subset of primes with positive \emph{relative} upper density:
$$ \limsup_{N \to \infty} \frac{|A \cap [1,N]|}{|P \cap [1,N]|} > 0.$$
Then $A$ contains infinitely many arithmetic progressions $n, n+r, n+2r$ of length three.
\end{theorem}

A key observation made in that paper was that one did not need very deep number-theoretic information about the structure of $A$ or $P$
to prove this result.  In fact, the same result holds not just for relatively dense subsets of primes, but relatively dense
subsets of \emph{almost primes} (numbers containing no small prime factors); we shall return to this point later.

In 2004, Ben Green and the author \cite{gt-primes} were able to extend this theorem to arbitrarily long progressions, by replacing
Fourier-analytic ideas with ergodic theory ones:

\begin{theorem}\cite{green-tao}
Let $A \subset P$ be a subset of primes with positive \emph{relative} upper density:
$$ \limsup_{N \to \infty} \frac{|A \cap [1,N]|}{|P \cap [1,N]|} > 0,$$
and let $k \geq 3$.  Then $A$ contains infinitely many arithmetic progressions $n, n+r, \ldots, n+(k-1)r$ of length $k$.  In particular,
the primes contain arbitrarily long arithmetic progressions.
\end{theorem}

At the time of writing, we are not able to obtain van der Corput's more precise asymptotic estimate on the \emph{number} of prime progressions of
arbitrary length $k$, but we are able to do so in the $k=4$ case; see Section \ref{unif-sec}.

In this expository article, we review briefly the methods of proof of Roth's theorem and Szemer\'edi's theorem for various values of $k$, focusing in 
particular on the cases $k=3$ and $k=4$ which are amenable to Fourier analysis and ``quadratic Fourier analysis'' respectively.  Then we discuss the recent extension of these theorems to the prime numbers.  There is substantial overlap between this survey and \cite{green-survey}.

\section{Progressions of length three}

We now discuss some proofs of Roth's theorem.  We first observe
that this theorem can be reformulated in one of two equivalent ``finitary'' settings: firstly as a statement about subsets of
long arithmetic progressions, and secondly as a statement about a large cyclic group.

We need some notation.  The interval $[a,b]$ shall always refer to the discrete interval $\{ n \in \Z: a \leq n \leq b \}$.  We use
$|A|$ to denote the cardinality of a finite set $A$.  If $A$ is a finite set and $f: A \to \C$ is a complex-valued function, we define
the \emph{expectation} $\E(f) = \E( f(n) |  n \in A )$ of $f$ to be the quantity
$$ \E( f(n) | n \in A ):= \frac{1}{|A|} \sum_{n\in A} f(n);$$
similarly, if $P(n)$ is a property pertaining to elements of $A$, we define the \emph{probability} of $P$ to be
$$ \P(A) = \P( P(n) | n \in A ) := \frac{1}{|A|} |\{ n \in A: P(n) \hbox{ is true} \}|,$$
and we define $1_{P}$ to be the indicator function of $P$, thus $1_P(n) = 1$ when $P(n)$ is true and $1_P(n) = 0$ otherwise.

\begin{theorem}[Roth's theorem, second version]  Let $0 < \delta \leq 1$.  Then there exists an $N_0 := N_0(\delta) > 1$ such that,
for any arithmetic progression $P \subset \Z$ of length at least $N_0$ and any subset $A \subset P$ of density $\P( n \in A: n \in P ) \geq \delta$, $A$ contains at least one arithmetic progression $n, n+r, n+2r$ of length three.
\end{theorem}

Note that the choice of progression $P$ is unimportant to this theorem; only the length is relevant.  This is because all progressions of a fixed length
are isomorphic to each other by an affine scaling map.  Thus one could set $P = [1,N]$ here for some $N \geq N_0$ with no loss of generality.

Henceforth let us call a function $f: A \to \C$ on a finite set $A$ \emph{bounded} if $|f(n)| \leq 1$ for all $n \in A$.

\begin{theorem}[Roth's theorem, third version]   Let $0 < \delta \leq 1$, and let $N \geq 1$ be a prime integer.  Let $f: \Z/N\Z \to [0,1]$ be a
non-negative bounded function with large mean
\begin{equation}\label{mean}
 \E( f(n) | n \in \Z/N\Z ) \geq \delta.
\end{equation}
Then we have
\begin{equation}\label{peon}
 \E( f(n) f(n+r) f(n+2r) | n, r \in \Z/N\Z ) \geq c(3,\delta) - o_{\delta}(1)
\end{equation}
for some $c(3,\delta) > 0$ depending only on $\delta$, where $o_{\delta}(1)$ is a quantity that depends on $\delta$ and $N$, and
for each fixed $\delta$ tends to zero as $N$ goes to infinity.
\end{theorem}

Before we prove any of these versions, let us first sketch why they are equivalent. 

\begin{proof}[Second version implies first version]  Let $A$ be a set of positive upper density.  Then there exists a $\delta > 0$ such that
$|A \cap [1,N]| \geq 2\delta N$ for infinitely many $N$.  Using this, one can find infinitely many disjoint intervals $[a_j, b_j]$ of length
$b_j - a_j \geq N_0(\delta)$ such that $A$ has density at least $\delta$ on these intervals:
$$ \P( n \in A: n \in [a_j,b_j] ) \geq \delta.$$
Applying the second version of Roth's theorem to each such interval we thus see $A$ has infinitely many progressions of length 3 as desired.
\end{proof}

\begin{proof}[First version implies second version] Suppose for contradiction that the second version failed.  Then we could find a 
$\delta > 0$ and sets $A_j \subset [1, N_j]$ (with $N_j \to \infty$) with $\P( n \in A_j: n \in [1,N_j] ) \geq \delta$ and with
each $A_j$ containing no arithmetic progressions of length 3.  By refining the sequence if necessary we may assume that the $N_j$
are increasing in $j$ (indeed we could make this sequence grow incredibly fast if desired).  If one then considers the 
set $A := \bigcup_{j=1}^\infty 2N_j + A_j$, then it is easy to show that $A$ has positive upper density but contains no arithmetic progressions,
a contradiction.
\end{proof}

\begin{proof}[Third version implies second version]  Let $p$ be a prime between $2N$ and $4N$ (which always exists by Bertrand's postulate).  
Let $\pi: [1,N] \to \Z/p\Z$ be the canonical injection of $[1,N]$ into $\Z/p\Z$.  If $A \subset [1,N]$ has density
$\P( n \in A: n \in [1,N] ) \geq \delta$, then the function $f := 1_{\pi(A)}$ on $\Z/p\Z$ is non-negative, bounded, and obeys the estimate
$$  \E( f(n) | n \in \Z/p\Z ) = \frac{|A|}{p} \geq \frac{|A|}{4N} \geq \delta/4.$$
Thus by the third version of Roth's theorem we have
$$ \E( f(n) f(n+r) f(n+2r) | n, r \in \Z/p\Z ) \geq c(3,\delta/4) - o_{\delta}(1).$$
Note that $f(n) f(n+r) f(n+2r)$ is non-zero only when $n = \pi(n')$, $n+r = \pi(n'+r')$, $n+2r = \pi(n'+2r')$ and $n' \in [1,N]$, $-N < r < N$, in 
which case this quantity is equal to 1.  Thus we have
$$ | \{ (n',r'): n',n'+r',n'+2r' \in A; n' \in [1,N]; -N < r' < N \} | \geq c(3,\delta/4) p^2 - o_{\delta}(p^2).$$
We can discard the $r'=0$ terms as they contribute $O(N) = o(p^2)$.  By symmetry we can then reduce to the
positive $r'$.  We thus have
$$ | \{ (n',r'): n',n'+r',n'+2r' \in A; n' \in [1,N]; 0 < r' < N \} | \geq c(3,\delta/4) p^2/2 - o_{\delta}(p^2).$$
If $N$ (and hence $p$) is sufficiently large, then the right-hand side is non-zero, and we have demonstrated the existence of a non-trivial
arithmetic progression of length three in $A$.  (In fact we have demonstrated $\geq c'(3,\delta) N^2$ such progressions for some
$c'(3,\delta) > 0$).
\end{proof}

\begin{proof}[Second version implies third version]  This argument is due to Varnavides \cite{varnavides}.  We first observe that to prove the
theorem, it suffices to do so when $f$ is a characteristic function $f = 1_A$.  This is because if $f$ is non-negative, bounded and obeys \eqref{mean}
then the set $A := \{ n \in \Z/N\Z: f(n) \geq \delta/2 \}$ must have density at least $\P( n \in A: n \in \Z/N\Z ) \geq \delta/2 )$.
Since we have the pointwise bound\footnote{This is somewhat crude.  A slightly better argument would be to select $A$ randomly, with each element $n \in \Z/N\Z$ having a probability of $f(n)$ to lie in $A$, and then take averages, but in practice this does not yield significantly better constants at the end.} from below $f \geq \frac{\delta}{2} 1_A$, we have
$$ \E( f(n) f(n+r) f(n+2r) | n, r \in \Z/N\Z ) \geq \frac{\delta^3}{8} \E( 1_A(n) 1_A(n+r) 1_A(n+2r) | n, r \in \Z/N\Z )$$
and so \eqref{peon} for $f$ would follow from \eqref{peon} for $A$ (with a slightly worse value of $c(3,\delta)$, namely
$\frac{\delta^3}{8} c(3,\delta/2)$).

It remains to verify \eqref{peon} for characteristic functions.
Let $M = M(\delta)$ be a large
integer depending on $\delta$ to be chosen later.  To prove \eqref{peon} it suffices to do so in the case $N \gg M$, since the case $N = O(M)$
is vacuous.  

The idea is to cover $\Z/N\Z$ uniformly by progressions $P_{ab} := \{ a + b, a+2b, \ldots, a+Mb \}$ of length $M$, where we allow $b$ to be zero.  Indeed we observe that for every $n \in \Z/N\Z$ there are exactly $NM$ pairs $(a,b) \in \Z/N\Z \times \Z/N\Z$ such that $n \in a + [1,M] \cdot b$ (this is easiest to see by choosing $b$ first).  Thus
\begin{align*}
\delta &\leq \P( n \in A | n \in \Z/N\Z )\\
&= \P( n \in A | n \in P_{ab}; (a,b) \in \Z/N\Z \times \Z/N\Z )\\
&= \E( \P(n \in A | n \in P_{ab} | (a,b) \in \Z/N\Z \times \Z/N\Z ).
\end{align*}
In particular, if we let $\Omega \subseteq \Z/N\Z \times \Z/N\Z$ be the set of pairs $(a,b)$ such that
$\P(n \in A | n \in P_{ab}) \geq \delta/2$, then we have
\begin{equation}\label{omega-dens}
 \P( (a,b) \in \Omega | (a,b) \in \Z/N\Z \times \Z/N\Z ) \geq \delta/2.
 \end{equation}
Now choose $M := N_0(\delta/2)$.  From the definition of $\Omega$ and the second form of Roth's theorem,
we see that for every $(a,b) \in \Omega$, the set $A \cap P_{ab}$
contains at least one non-trivial arithmetic progression $n, n+r, n+2r$ of length three.  In particular
we have
$$ \P( n, n+r, n+2r \in A | n, n+r, n+2r \in P_{ab}; r \neq 0 ) \geq M^{-2}$$
since the number of progressions $n, n+r, n+2r$ in $P_{ab}$ is at most $M^2$.

Now observe that every progression $n, n+r, n+2r \in \Z/N\Z$ with $r \neq 0$ is contained in exactly the same number of progressions $P_{ab}$, since
they are all isomorphic using affine scaling maps (here we use that $N = |\Z/N\Z|$ is prime).  Thus we have
$$ \P( n, n+r, n+2r \in A | n, n+r, n+2r \in \Z/N\Z; r \neq 0 ) \geq M^{-2}$$
In particular (adding in the $r=0$ case) we have
$$ \E( \prod_{j=0}^2 1_A(x+jr) | x \in \Z/N\Z; r \in \Z/N\Z )
\geq  M^{-2} - o(1)$$
which gives \eqref{peon} as desired (with $c(3,\delta) = M^{-2} = N_0(\delta/2)^{-2}$ for characteristic
functions, and hence $c(3,\delta) = \frac{\delta^3}{8} N_0(\delta/4)^{-2}$ for arbitrary functions).
\end{proof}

In light of these equivalent formulations, it is natural to introduce the Lebesgue spaces $L^p(\Z/N\Z)$ for $1 \leq p \leq \infty$, 
defined as the complex-valued functions on $\Z/N\Z$ equipped with the norm
$$ \| f \|_{L^p(\Z/N\Z)} := \E( |f|^p )^{1/p} = (\frac{1}{N} \sum_{n \in \Z/N\Z} |f(n)|^p)^{1/p}$$
and to introduce the trilinear form $\Lambda_3: L^{p_1}(\Z/N\Z) \times L^{p_2}(\Z/N\Z) \times L^{p_3}(\Z/N\Z) \to \C$ by
\begin{equation}\label{l3-def}
\Lambda_3(f,g,h) := \E( f(n) g(n+r) h(n+2r) | n,r \in \Z/N\Z ).
\end{equation}
Here we always assume $N$ to be a large prime (in particular, it is odd).
Thus the third version of Roth's theorem can be reformulated as follows: if $f \in L^\infty(\Z/N\Z)$ is a non-negative function obeying
the bounds
$$ 0 < \delta \leq \|f\|_{L^1(\Z/N\Z)} \leq \|f\|_{L^{\infty}(\Z/N\Z)} \leq 1$$
then
\begin{equation}\label{t3}
 \Lambda_3(f,f,f) \geq c(3,\delta) - o_{\delta}(1)
 \end{equation}
for some $c(3,\delta) > 0$.
Note that the task here is to obtain \emph{lower bounds} on the form $\Lambda_3(f,f,f)$ rather than upper bounds, which are considerably easier to
obtain.  For instance, from multilinear interpolation (or Young's inequality) it is easy to establish the upper bounds
\begin{equation}\label{young}
 |\Lambda_3(f,g,h)| \leq \| f \|_{L^p} \|g\|_{L^q} \|h\|_{L^r}
\end{equation}
whenever $1 \leq p,q,r \leq \infty$ and $\frac{1}{p} + \frac{1}{q} + \frac{1}{r} \leq 2$; here $f,g,h$ are arbitrary
complex-valued functions.  Note that the non-negativity of $f$ and of $\Lambda_3$ (i.e. $\Lambda_3(f,g,h)$ is non-negative whenever $f,g,h$
are non-negative) is crucial, since without this one could not even obtain the trivial bound\footnote{There is also the slightly better
trivial bound $\Lambda_3(f,f,f) \geq \| f\|_{L^3(\Z/N\Z)}^3/N$ coming from the $r=0$ term in \eqref{l3-def}, but this lower bound is $o(1)$ and is thus not significantly better than the trivial bound of 0.}
$\Lambda_3(f,f,f) \geq 0$, let alone \eqref{t3}.

At first glance it does not appear that upper bounds such as \eqref{young} are useful for proving lower bounds of the type \eqref{t3}.  
However, one can use the multilinearity of $\Lambda_3$ to convert upper bounds to lower bounds as follows.  Without loss of generality
we may take $\E(f)$ to be equal to $\delta$ (since if $\E(f) >\delta$ we may simply decrease $f$ and hence $\Lambda_3(f,f,f)$.
We decompose\footnote{This is of course a very simple decomposition.  Later on we shall use more sophisticated decompositions, which
can be viewed as ``arithmetic'' versions of the Calder\'on-Zygmund decomposition in harmonic analysis.} $f$ into 
a ``good function'' $g := \E(f) = \delta$ and a ``bad function'' $b := f - \E(f)$, and then we can split $\Lambda_3(f,f,f)$ into eight components:
$$ \Lambda_3(f,f,f) = \Lambda_3(g,g,g) + \ldots + \Lambda_3(b,b,b).$$
The first term can be computed explicitly, and can be viewed as a main term:
$$ \Lambda_3(g,g,g) = \Lambda_3(\delta,\delta,\delta) = \delta^3.$$
Thus if one can obtain \emph{upper} bounds on the magnitude of the remaining seven terms which add up to less than $\delta^3$, then one can hope to
prove \eqref{t3}.  The bound \eqref{young} turns out to be too weak to do this, unless $\delta$ is very close to 1 (e.g. if $\delta > 2/3$);
however, one can do better by replacing the Lebesgue norms with some additional norms, based on the \emph{Fourier transform}
$$ \hat f(\xi) := \E( f(x) e_N(-x \xi) | x \in \Z/N\Z ),$$
where $e_N: \Z/N\Z \to S^1$ is the character $e_N(x) := \exp(2\pi i x/N)$.
From the Fourier inversion formula
$$ f(x) = \sum_{\xi \in \Z/N\Z} \hat f(\xi) e_N(x\xi)$$
we see that
$$
\Lambda_3(f,g,h) = \sum_{\xi_1, \xi_2, \xi_3 \in \Z/p\Z} \hat f(\xi_1) \hat g(\xi_2) \hat f(\xi_3)
\E( e_N( n\xi_1 + (n+r) \xi_2 + (n+2r) \xi_3 ) | n,r \in \Z/p\Z ).$$
The expectation on the right-hand side equals 1 when $\xi_1 = \xi_3$ and $\xi_2 = -2 \xi_1$, and
equal to zero otherwise.  Thus we have the identity
$$ 
\Lambda_3(f,g,h) = \sum_{\xi \in \Z/p\Z} \hat f(\xi) \hat g(-2\xi) \hat h(\xi).$$
From the Plancherel identity
$$ \| f \|_{L^2(\Z/N\Z)} = \| \hat f \|_{l^2(\Z/N\Z)}$$
and H\"older's inequality, we thus have the estimate
\begin{equation}\label{lambda3-est}
|\Lambda_3(f,g,h)| \leq \| f \|_{L^2(\Z/N\Z)} \|g\|_{L^2(\Z/N\Z)} \| \hat h \|_{l^\infty(\Z/N\Z)}
\end{equation}
and similarly for permutations.  We also have the variant
\begin{equation}\label{lambda4-est}
|\Lambda_3(f,g,h)| \leq \| f \|_{L^2(\Z/N\Z)} \|\hat g\|_{l^4(\Z/N\Z)} \| \hat h \|_{l^4(\Z/N\Z)}
\end{equation}

This leads to the following criterion to ensure $\Lambda_3(f,f,f)$ is positive.

\begin{proposition}\label{gbgb}  Let $f \in L^\infty(\Z/N\Z)$ have a decomposition of the form $f = g + b$, where
\begin{equation}\label{g-bound}
 \| g \|_{L^\infty(\Z/N\Z)}, \| b \|_{L^\infty(\Z/N\Z)} = O(1); \quad
\| g \|_{L^\infty(\Z/N\Z)}, \| b \|_{L^\infty(\Z/N\Z)} = O(\delta).
\end{equation}
Then we have the estimates
\begin{equation}\label{l3-infty}
\Lambda_3(f,f,f) = \Lambda_3(g,g,g) + O( \delta \| \hat b \|_{l^\infty(\Z/N\Z)} )
\end{equation}
and
$$ \Lambda_3(f,f,f) = \Lambda_3(g,g,g) + O( \delta^{5/4} \| \hat b \|_{l^4(\Z/N\Z)} ).$$
\end{proposition}

\begin{remark} Interestingly, estimates of this type (after being suitably localized in phase space) have proven to be crucial in
recent progress in understanding the bilinear Hilbert transform (see e.g. \cite{laceyt1}), or at least in understanding the contribution of individual ``trees''
to that transform.  Indeed there is some formal similarity between the trilinear form $\Lambda_3$ and the trilinear form $\Lambda(f,g,h) := p.v. \int\int f(x+t)g(x-t)h(x) \frac{dx dt}{t}$ associated to the bilinear Hilbert transform.
\end{remark}

\begin{proof}  From the hypotheses we have
$$ \| g \|_{L^2(\Z/N\Z)}, \| b \|_{L^2(\Z/N\Z)} = O(\delta^{1/2})$$
and hence by Plancherel
$$ \| \hat g \|_{l^2(\Z/N\Z)}, \| \hat b \|_{l^2(\Z/N\Z)} = O(\delta^{1/2}).$$
On the other hand, from the $L^1$ bounds on $g$ and $b$ we have
$$ \| \hat g \|_{l^\infty(\Z/N\Z)}, \| \hat b \|_{l^\infty(\Z/N\Z)} = O(\delta)$$
and so by H\"older's inequality
$$ \| \hat g \|_{l^4(\Z/N\Z)}, \| \hat b \|_{l^4(\Z/N\Z)} = O(\delta^{3/4}).$$
The claims now follow by decomposing $\Lambda_3(f,f,f)$ into eight pieces as before, setting aside $\Lambda_3(g,g,g)$ as a main term,
and using \eqref{lambda3-est}, \eqref{lambda4-est}
(and permutations thereof) to estimate all the remaining pieces (which involve at least one copy of $b$).
\end{proof}

This suggests the following strategy: in order to obtain a non-trivial lower bound on $\Lambda_3(f,f,f)$, we should obtain a splitting
$f=g+b$ obeying the bounds \eqref{g-bound}
where the ``good'' function $g$ already has a large value of $\Lambda_3(g,g,g)$ (thus we shall presumably want $g$ to be non-negative), and the ``bad'' function $b$
has a small Fourier transform, either in $l^\infty$ norm or $l^4$ norm.  Note that up to polynomial factors of $\delta$, the two norms are somewhat
equivalent, as one can easily establish the estimates
\begin{equation}\label{lest}
 \| \hat b \|_{l^\infty} \leq \| \hat b \|_{l^4} \leq \| \hat b \|_{l^\infty}^{1/2} \| \hat b \|_{l^2}^{1/2}
\leq C \delta^{1/4} \| \hat b \|_{l^\infty}^{1/2}.\end{equation}
In the original arguments involving Roth's theorem, the $l^\infty$ norm on the Fourier coefficients was used, but as we shall see later, it
is the $l^4$ norm which is easier to generalize to ``higher order'' Fourier analysis, which will be necessary to treat the $k \geq 4$
case.  Let us rather informally call a function $b$ which obeys bounds such as \eqref{g-bound} 
\emph{linearly uniform} if the Fourier transform $\hat b$ is very small in either $l^\infty$ or $l^4$; we see from \eqref{lest} that
 it is not terribly important which norm we choose here.  The reason for this terminology is that a linearly uniform function $b$ is one which is
 uniformly distributed with respect to linear phase functions $e_N(x \xi)$, in the sense that the inner product of $b$ with such
 functions is small.  (This rather vague statement can be made more precise using Weyl's criterion for uniform distribution).

We have already indicated one such candidate for a decomposition, namely the decomposition $f = g + b$ into the expectation $g := \E(f)$ and
the expectation-free $b := f - \E(f)$ components of $f$.  Certainly this decomposition obeys the bounds \eqref{g-bound}, and the
value of $\Lambda_3(g,g,g) \geq \delta^3$ is moderately large.  However, at this stage we do not have very good bounds on $\|\hat b \|_{l^\infty(\Z/N\Z)}$ or $\| \hat b \|_{l^4(\Z/N\Z)}$; the best bounds we have on these quantities are $O(\delta)$ and $O(\delta^{3/4})$
respectively, and thus the error term can dominate the main term.  (Indeed, there certainly exist functions $f$ for which $\Lambda_3(f,f,f)$ is
significantly different from $\Lambda_3(g,g,g)$; consider for instance $f = 1_{[1,\delta N]}$, in which the former quantity is comparable to $\delta^2$ and the latter is comparable to $\delta^3$).

However, we can at least eliminate one case, in which $b = f - \E(f)$ is sufficiently linearly uniform (for instance if $\| \hat b \|_\infty \leq \delta^2/100$).  The question is then what to do in the remaining cases, when $b = f - \E(f)$ is not sufficiently linear uniform.  The strategy is then
to \emph{convert} the lack of linear uniformity from a liability to an asset, by showing that this lack of uniformity implies some additional structure which one can exploit to improve the situation.  The known proofs of Roth's theorem (or more generally Szemer\'edi's theorem) differ on exactly what this additional structure could be, and how to exploit it, but they essentially fall into one of two categories\footnote{Szemer\'edi's proof of
Szemer\'edi's theorem in \cite{szemeredi} is a blend of the density increment and energy increment arguments.}:

\begin{itemize}
\item A \emph{density increment argument} seeks to use the lack of uniformity in $b$ to pass from $\Z/N\Z$ (or $[1,N]$) to a smaller object on which
the function $f$ (or the set $A$) has a larger density.  One then iterates this procedure until uniformity is obtained; this algorithm terminates since
the density is bounded.

\item An \emph{energy increment argument} seeks to use the lack of uniformity in $b$ to improve the decomposition $f=g+b$, replacing the good function $g$ by a function of larger energy ($L^2$ norm).  One then iterates this procedure until uniformity is obtained; this algorithm terminates
since the energy is bounded.
\end{itemize}

Both approaches are important to the theory, as they have different strengths and weaknesses.  We illustrate this by giving two proofs of
Roth's theorem, one for each of the above approaches.  But we shall need some additional notation first; this notation may seem somewhat
cumbersome for this application, but will become very convenient when we discuss the case of larger $k$ in later sections.

\begin{definition}[$\sigma$-algebras] Let $X$ be a finite set (such as $\Z/N\Z$ or $[1,N]$).  A \emph{$\sigma$-algebra} $\B$ in $X$ is any collection
of subsets of $X$ which contains the empty set $\emptyset$ and the full set $X$, and is closed under
complementation, unions and intersections.  We define the \emph{atoms} of a $\sigma$-algebra to be the minimal non-empty elements of $\B$ (with respect to set inclusion); it is clear that the atoms in $\B$ form a partition of $X$, and $\B$ consists precisely of arbitrary unions of its atoms (including the empty union $\emptyset$); thus there is a one-to-one correspondence between $\sigma$-algebras and partitions of $X$.
A function $f: X \to \C$ is said to be \emph{measurable} with respect to a $\sigma$-algebra $\B$ if all the level 
sets of $f$ lie in $\B$, or equivalently if $f$ is constant on each of the atoms of $\B$.  
We define $L^2(\B)$ be the space of $\B$-measurable functions, equipped with the Hilbert space inner product $\langle f, g \rangle_{L^2(X)} := \E( f \overline{g} )$.  We can then define the conditional expectation operator $f \mapsto \E(f|\B)$ to
be the orthogonal projection of $L^2(X)$ to $L^2(\B)$.  An equivalent definition of conditional expectation is
$$ \E(f|\B)(x) := \E( f(y) | y \in \B(x) )$$
for all $x \in X$, where $\B(x)$ is the unique atom in $\B$ which contains $x$.  
It is clear that conditional
expectation is a linear self-adjoint orthogonal projection on $L^2(\mathbb{Z}_N)$, preserves non-negativity, expectation, and constant functions.  In particular it maps bounded functions to bounded functions.  If $\E(f|\B)$ is zero we say
that $f$ is \emph{orthogonal to} $\B$.

If $\B$, $\B'$ are two $\sigma$-algebras, we use $\B \vee \B'$ to denote the $\sigma$-algebra generated by $\B$ and $\B'$ (i.e. the $\sigma$-algebra whose atoms are the intersections of atoms in $\B$ with atoms in $\B'$).  
\end{definition}

\begin{proof}[Density increment proof of Roth's theorem]  We now give what is essentially Roth's original argument, though not using Roth's original language (in particular, we give the sigma algebras of Bohr sets significantly more prominence in the argument).

It is more convenient to work with the second formulation of Roth's theorem.  Let $\delta > 0$, and let $N$ be a sufficiently large number
depending on $\delta$.  Let $P_0$ be a progression of length $N$, and let $A$ be a subset of $P_0$ of density at least $\delta$.  Our task is to prove that $A$ contains at least one
arithmetic progression.

Without loss of generality we may take $P_0 = [1,N]$.
Set $\delta_0 := \P( n \in A: n \in [1,N])$, thus $\delta \leq \delta_0 \leq 1$.

Choose a prime $p$ between $2N$ and $4N$.  We embed $[1,N]$ into $\Z/p\Z$ in the obvious manner, thus identifying $A$ with a subset of $\Z/p\Z$,
of density at least $\delta/4$.  Let us let $f: \Z/p\Z \to \R$ be defined by setting $f(x) := 1_A(x)$ when $x \in [1,N]$ and
$f(x) = \delta_0$ otherwise; observe that $\E(f) = \delta_0$ by construction.  We then split $f = g+b$, where $g := \E(f) \geq \delta_0$ and $b := f-\E(f) = 1_A - \delta_0 1_{[1,N]}$.

There are two cases, depending on whether $b$ is linearly uniform or not.  Suppose first that $b$ is linearly uniform in the sense
that $\|\hat b\|_{l^\infty(\Z/p\Z)} \leq c \delta^2$ for some small absolute constant $0 < c \ll 1$; this is the ``easy case''.  Since $\Lambda_3(g,g,g) = \E(f)^3 \geq \delta_0^3 \geq \delta^3$, we see from \eqref{l3-infty} that $\Lambda_3(f,f,f) \geq c' \delta^2$ for some absolute constant $c' > 0$ (if $c$ is chosen
sufficiently small).  By definition of $f$ and $\Lambda_3$, this means that
$$ \P( n, n+r, n+2r \in A | (n, r) \in \Z/p\Z \times \Z/p\Z ) \geq c' \delta^2.$$
The contribution of the $r=0$ case is at most $O(1/p) = O(1/N)$.  Thus if $N$ is large enough, we thus see that there exists at least one
pair $(n,r) \in \Z/p\Z \times \Z/p\Z$ with $r \neq 0$ such that $n,n+r,n+2r$ in $A$.  Since $A \subseteq [1,N]$, this forces $n \in [1,N]$
and $1 \leq |r| \leq N$.  Since $p > 2N$, this implies that $A$ (thought now as a subset of
$\Z$ rather than $\Z/p\Z$) also contains a non-trivial arithmetic progression $n, n+r, n+2r$, as claimed.

Now suppose we are in the ``hard case'' where $b$ is not linearly uniform, then there exists a frequency $\xi \in \Z/p\Z$ such that $|\hat b(\xi)| \geq c\delta^2$.  By definition
of $b$ and the Fourier transform, we thus have
$$ |\E( (1_A(n) - \delta_0 1_{[1,N]}(n)) e_p(-n\xi) | n \in \Z/p\Z )| \geq c \delta^2.$$
Transferring this back from $\Z/p\Z$ to $[1,N]$, we obtain
$$ |\E( (1_A(n) - \delta_0) e_p(-n\xi) | n \in [1,N] )| \geq c \delta^2$$
(with a slightly different constant $c$).  If we let $\chi: [1,N] \to \C$ be the linear phase function
$\chi(n) := e_p(n\xi)$, we see that $1_A-\delta_0$ thus has some correlation with $\chi$:
\begin{equation}\label{adc}
|\langle 1_A-\delta_0, \chi \rangle_{L^2([1,N])}| \geq c \delta^2.
\end{equation}

Now let $0 < \eps \ll 1$ be a small quantity depending on $\delta$ to be chosen later.  We partition the complex plane $\C = \bigcup_{Q \in \Q_\eps} Q$
into squares of side-length $\eps$ in the standard manner (i.e. the corners of the square lie in the lattice $\eps \Z^2$), and let $\B_{\eps,\chi}$ be the $\sigma$-algebra on $[1,N]$ 
generated by the atoms $\{ \chi^{-1}(Q): Q \in \Q_\eps \}$; sets of this type are also known as \emph{Bohr sets}.
Observe that there are only $O(1/\eps)$ non-empty atoms.  
Then on each atom, $\chi$ can only vary by at most $O(\eps)$, and thus we have
the pointwise estimate
$$ \chi - \E(\chi | \B_{\eps,\chi}) = O(\eps).$$
Since $1_A - \delta_0$ is bounded, we thus see from \eqref{adc} and the triangle inequality that
$$
|\langle 1_A-\delta_0, \E(\chi|\B_{\eps,\chi}) \rangle_{L^2([1,N])}| \geq c \delta^2 - O(\eps).
$$
Since conditional expectation is self-adjoint, we have
$$ \langle 1_A-\delta_0, \E(\chi|\B_{\eps,\chi}) \rangle_{L^2([1,N])} =
\langle \E(1_A-\delta_0|\B_{\eps,\chi}), \chi \rangle_{L^2([1,N])},$$
and thus by boundedness of $\chi$
$$ \| \E(1_A - \delta_0|\B_{\eps,\chi}) \|_{L^1([1,N])} \geq c \delta^2 - O(\eps).$$
If we choose $\eps := c' \delta^2$ for some suitably small absolute constant $0 < c' \ll 1$, the left-hand side is at least $c\delta^2/2$.
Now observe that $\E(1_A - \delta_0|\B_{\eps,\chi})$ has mean zero:
$$ \E( \E(1_A - \delta_0|\B_{\eps,\chi}) ) = \E( 1_A - \delta_0 ) = \E(1_A) - \delta_0 = \delta_0 - \delta_0 = 0.$$
Thus we see that the positive part of $\E(1_A - \delta_0|\B_{\eps,\chi})$ is large:
$$ \E( \E(1_A - \delta_0|\B_{\eps,\chi})_+ ) \geq c \delta^2/4.$$
Now recall that $\B_{\eps,\chi}$ is generated by $O(1/\eps) = O(\delta^{-2})$ non-empty atoms.  By definition of conditional expectation
and the pigeonhole principle, we can thus find some atom $\chi^{-1}(Q)$ of $\B_{\eps,\chi}$ of density at least $c'' \delta^4$
such that $1_A - \delta_0$ is biased on this atom:
$$ \E( 1_A(n) - \delta_0 | n \in \chi^{-1}(Q) ) \geq c''' \delta^2,$$
and thus
\begin{equation}\label{pnaq}
 \P( n \in A | n \in \chi^{-1}(Q) ) \geq \delta_0 + c''' \delta^2.
 \end{equation}
This is a \emph{density increment}; $A$ is denser on $\chi^{-1}(Q)$ than it is on $[1,N]$.  However, $\chi^{-1}(Q)$ is a Bohr set instead
of an arithmetic progression.  However, the Bohr set is in some sense ``very close'' to an arithmetic progression in the sense that it can be
covered quite efficiently by somewhat long arithmetic progressions\footnote{This step is not particularly efficient when it comes to quantitative constants.  A more
refined argument of Bourgain \cite{bourgain-triples} works entirely with Bohr sets rather than arithmetic progressions, and obtains the best bounds
on $N_0(\delta)$ to date (namely $N_0(\delta) \leq C \delta^{-C/\delta^2}$).}.  This can be seen as follows.  By the pigeonhole principle, one
can find an integer $1 \leq q \leq \sqrt{N}$ such that
$$ \| q \frac{\xi}{p} \| \leq \frac{1}{\sqrt{N}},$$
where $\|x\|$ denotes the distance of $x$ to the nearest integer.  From this one easily observes that if $n \in \chi^{-1}(Q)$, then there is
an arithmetic progression containing $n$ of spacing $q$ and length comparable to $\eps \sqrt{N}$ which is completely contained in $\chi^{-1}(Q)$.
In particular, one can partition $\chi^{-1}(Q)$ into disjoint arithmetic progressions, each of length comparable to $\eps \sqrt{N} \geq C^{-1} \delta^2 \sqrt{N}$.
From \eqref{pnaq} and the pigeonhole principle, we thus see that at least one of these progressions $P_1$ has large density:
$$ \P( n \in A | n \in P_1 ) \geq \delta_0 + c''' \delta^2.$$

To summarize, we had started with a subset $A$ of a progression $P_0$ of length $N$ which had density $\delta_0$, and concluded that either $A$ contained an arithmetic progression, or there was a sub-progression $P'$ of length at least $C^{-1} \delta^2 \sqrt{N}$ where $A$ has density
$\delta_1 > \delta_0 + c''' \delta^2$ for some absolute constant $c''' > 0$.  We can then pass to this progression $P'$ and repeat the argument
(note that we can make $C^{-1} \delta^2 \sqrt{N}$ as large as we please by requiring $N$ to be sufficiently large).  The density can only increase by $c''' \delta^2$ by at most $O(1/\delta^2)$ times\footnote{One can improve this to $O(1/\delta)$ by observing that the density increment of $c''' \delta^2$ can be refined to $c''' \delta_0^2$.}, and so this argument must eventually yield a non-trivial arithmetic of length three in $A$.
\end{proof}

\begin{proof}[Energy increment proof of Roth's theorem]  We now give an energy increment proof of Roth's theorem, inspired by
arguments of Furstenberg \cite{furst}, Bourgain \cite{bourg2}, and Green \cite{green}, as well as later arguments by Green and the author in \cite{green-tao}, \cite{tao:szemeredi}.  This is not the shortest such proof, nor the most efficient as far as explicit bounds are concerned, but it is a proof which has a relatively small reliance on Fourier analysis and thus which generalizes fairly easily to general $k$.  The structure of this argument, and the concepts introduced, are particularly crucial when establishing long arithmetic progressions in the primes.

We shall use the third formulation of Roth's theorem; unlike the preceding proof, we will not oscillate back and forth between progressions
and cyclic groups, but remain in a fixed cyclic group $\Z/N\Z$ throughout.  
Thus, we let $N$ be a large prime, and let $f$ be a bounded non-negative function on $\Z/N\Z$
obeying the bound \eqref{mean}.  Our task is to prove \eqref{t3}.

We need some additional notation.  

\begin{definition}[Almost periodic functions]\label{ap-def}  A \emph{linear phase function} is a function $\chi: \Z/N\Z \to \C$ of the form $\chi(n) = e_N(n\xi)$ for 
some $\xi \in \Z/N\Z$, which we refer to as the \emph{frequency} of $\chi$.  If $K > 0$, then an \emph{$K$-quasiperiodic function} is a function $f$ of the form $\sum_{j=1}^K c_j \chi_j$, where each $\chi_j$ is a linear phase function (not necessarily distinct), and $c_j$ are scalars such that $|c_j| \leq 1$.  If $\sigma > 0$, then an \emph{$(\sigma,K)$-almost periodic function} is a function $f$ such that $\|f-f_{QP}\|_{L^2(\Z/N\Z)} \leq \sigma$ for some $K$-quasiperiodic function $f_{QP}$.
\end{definition}

Observe that if $f$ and $g$ are $(\sigma, K)$-almost periodic functions, then $fg$ is a $(2\sigma, K^2)$-almost periodic function (taking $(fg)_{K^2} := f_K g_K$).  

A key property of almost periodic functions is that one can obtain non-trivial lower bounds on the $\Lambda_3$ quantity:

\begin{lemma}[Almost periodic functions are recurrent]\label{ap-recur}  Let $0 < \delta < 1$ and $0 < \sigma \leq \delta^3/100$, and $f$ be an bounded non-negative $(\sigma,K)$-almost periodic function obeying \eqref{mean}.  Then we have
$$ \Lambda_3(f,f,f) \geq c(K,\delta) - o_{n,\delta}(1)$$
for some $c(K,M,\delta) > 0$ (the key point here being that this quantity is independent of $N$).
\end{lemma}

\begin{proof} 
Let $f_{QP} = \sum_{j=1}^K c_j \chi_j$ be the $K$-quasiperiodic function approximating $f$, and let $0 < \eps$ be a small
number (depending on $K$, $\delta$) to be chosen later.  Let $\xi_1,\ldots,\xi_K$ be the frequencies associated to the characters $\chi_1,\ldots,\chi_K$.  By Dirichlet's simultaneous approximation by rationals theorem (or the pigeonhole principle), we have
\begin{equation}\label{rx}
 \P( \| r \xi_j \| \leq \eps \hbox{ for all } 1 \leq j \leq K | r \in \Z/N\Z) \geq c(\eps,K) 
\end{equation}
for some $c(\eps,K)> 0$ independent of $N$.  Next, observe from the triangle inequality that if $r$ is as above, then
$$ \| T^r f_{AP} - f_{AP} \|_{L^2(\Z/N\Z)} \leq C(K) \eps$$
where $T^r$ is the shift map $T^r f(x) := f(x+r)$.  From this and the triangle inequality, we conclude
$$ \| T^r f - f \|_{L^2(\Z/N\Z)} \leq \delta^3 / 10 + C(K) \eps,$$
and by another application of $T^r$, we have
$$ \| T^{2r} f - T^r f \|_{L^2(\Z/N\Z)} \leq \delta^3 / 10 + C(K) \eps.$$
From this and the boundedness of $f$, we conclude that
$$  \| f T^r f T^{2r} f - f^3 \|_{L^1(\Z/N\Z)} \leq \delta^3/2 + C(K) \eps,$$
but from the bounded non-negativity of $f$, \eqref{mean}, and H\"older's inequality we have
$$ \| f^3 \|_{L^1(\Z/N\Z)} \geq \| f \|_{L^1(\Z/N\Z)}^3 \geq \delta^3$$
and hence (by positivity of $f$)
$$ \E( f T^r f T^{2r} f(n) | n \in \Z/N\Z ) \geq \delta^3/2 - C(K) \eps.$$
If we choose $\eps$ small enough depending on $\delta$ and $M$, we thus have
$$ \E( f T^r f T^{2r} f(n) | n \in \Z/N\Z ) \geq \delta^3/4.$$
Averaging over all $r$, using \eqref{rx} and the non-negativity of $f$, we obtain
$$ \E( f T^r f T^{2r} f(n) | n,r \in \Z/N\Z ) \geq \delta^3 c(\eps,K)/4.$$
But the left-hand side is nothing more than $\Lambda_3(f, f, f)$.  The claim follows.
\end{proof}
 
To exploit the above result we shall need to approximate a general function $f$ by an almost periodic function, plus a linearly
uniform error.  The first step in this strategy shall be to construct $\sigma$-algebras such that the measurable functions in this algebra
are all almost periodic.

\begin{lemma}  Let $0 < \eps \ll 1$ and let $\chi$ be a linear phase function.  Then there exists a $\sigma$-algebra $\B_{\eps,\chi}$ such that
$\| \chi - \E(\chi|\B_{\eps,\chi}) \|_{L^\infty} \leq C \eps$, and such that for every $\sigma>0$, there exists $K = K(\sigma,\eps)$
such that every function $f$ which is measurable with respect to $\B_{\eps,\chi}$ and obeys the bound
$\|f\|_{L^\infty(\Z/N\Z)} \leq 1$ is $(\sigma,K)$-almost periodic.
\end{lemma}

\begin{proof}  We use a random construction, constructing a $\sigma$-algebra which has the stated properties with non-zero probability.  
Let $\alpha$ be a randomly selected element of the unit square in the complex plane, and
let $\B_{\eps,\chi}$ be the
$\sigma$-algebra with atoms of the form $\{ \chi^{-1}(Q): Q \in \Q_\eps + \eps \alpha \}$.  Then as in the previous proof of Roth's theorem,
we have $\| \chi - \E(\chi|\B_{\eps,\chi}) \|_{L^\infty} \leq C \eps$.  Now we prove the approximation claim.  It suffices to verify the claim
for $\sigma = 2^{-n}$ for some integer $n \gg 1$, with probability $1 - O(\sigma)$.  Also, since $\B_{\eps,\chi}$ has at most $C(\eps)$ atoms, it suffices to verify the claim when $f$ is the indicator function of one of these atoms $A$, with probability $1 - O( C(\eps)^{-1} \sigma )$.

The function $f$ can be rewritten as $f(x) = 1_Q( \chi(x) - \eps \alpha )$.  We can use the Weierstrass approximation theorem to approximate $1_Q(z)$ on the disk $z = O(1/\eps)$ by a polynomial $P(z,\overline{z})$ involving at most $C(\sigma,\eps)$ terms
and with coefficients bounded by $C(\sigma,\eps)$ such that $|P|$ is bounded by $1$ in this disk, and 
$1_Q(z) - P(z,\overline{z}) = O( C^{-1} \sigma )$ for all $z$ in this
disk, except for a set of measure $O( C(\eps)^{-2} \sigma^2 )$.  A standard randomization argument then allows us to assert that
$$ \| 1_Q( \chi(x) - \eps \alpha ) - P( \chi(x) - \eps \alpha, \overline{\chi(x) - \eps \alpha} ) \|_{L^2(\Z/N\Z)} \leq \sigma$$
with probability $O( C(\eps)^{-1} \sigma)$.  But $P( \chi(x) - \eps \alpha, \overline{\chi(x) - \eps \alpha} )$ can be written
as the linear combination of at most $C(\eps,\sigma)$ characters, with coefficients at most $C(\eps,\sigma)$, and is thus $C(\eps,\sigma)$-quasiperiodic (one can reduce the coefficients to be less than 1 by repeating characters as necessary).  The claim 
follows.
\end{proof}

One can concatenate these $\sigma$-algebras together.  If $\B_1,\ldots,\B_n$ are $\sigma$-algebras, we let $\B_1 \vee \ldots \vee \B_n$ be the
smallest $\sigma$-algebra which contains all of them.

\begin{corollary}\label{sigma-cor}  Let $0 < \eps_1,\ldots,\eps_n \ll 1$ and let $\chi_1,\ldots,\chi_n$ be linear phases.  Let $\B_{\eps_1,\chi_1}, \ldots, \B_{\eps_1,\chi_n}$ be
the $\sigma$-algebras arising from the above corollary.  Then for every $\sigma>0$, there exists $K = K(n, \sigma,\eps_1,\ldots,\eps_n)$
such that every function $f$ which is measurable with respect to $\B_{\eps_1,\chi_1} \vee \ldots \vee \B_{\eps_1,\chi_n}$ and obeys the bound
$\|f\|_{L^\infty(\Z/N\Z)} \leq 1$ is $(\sigma,K)$-almost periodic.
\end{corollary}

\begin{proof} Since the number of atoms in this $\sigma$-algebra is at most $C(n, \eps_1,\ldots,\eps_n)$, it suffices to verify this when $f$ is the indicator
function of a single atom.  But then $f$ is the product of $n$ indicator functions from atoms in $\B_{\eps_1,\chi_1}, \ldots, \B_{\eps_n,\chi_n}$,
and the claim follows from the preceding lemma and the previously made observation that the product of almost periodic functions is almost periodic.
\end{proof}

The significance of these $\sigma$-algebras is not only that they contain functions which are almost periodic and hence have non-trivial bounds
on the $\Lambda_3$ form, but also that they capture ``obstructions to linear uniformity'':

\begin{lemma}[Non-uniformity implies structure]\label{sigma-grow}  Let $b$ be a bounded function such that $\| \hat b \|_{l^\infty} \geq \sigma > 0$, and let $0 < \eps \ll \sigma$.  Then
there exists a linear phase function $\chi$ with associated $\sigma$-algebra $\B_{\eps,\chi}$ such that
$$ \| \E( b | \B_{\eps,\chi} ) \|_{L^2(\Z/N\Z)} \geq C^{-1} \sigma.$$
\end{lemma}

This is proven by a repetition of the arguments used in the first proof of Roth's theorem, and we leave it to the reader.

We can now assemble all these ingredients together to prove Roth's theorem.  The major step here is a \emph{structure theorem} which decomposes an arbitrary function into an almost periodic piece and a linearly uniform piece.

\begin{proposition}[Quantitative Koopman-von Neumann theorem]\label{qkn}  Let $F: \R^+ \times \R^+ \to \R^+$ be an arbitrary function, let $0 < \delta \leq 1$, and let $f$ be any bounded non-negative
function on $\Z/N\Z$ obeying \eqref{mean}.  Let $\sigma := \delta^3/100$.  Then there exists a quantity $0 < K \leq C(F, \delta)$ and a decomposition $f = g+b$, where $g$ is bounded, non-negative, has mean $\E(g) = \E(f)$, and $(\sigma,K)$-almost periodic, and $b$ obeys the bound 
\begin{equation}\label{bf}
\| \hat b \|_{l^\infty} \leq F(\delta,K).
\end{equation}
\end{proposition}

\begin{proof}  We apply the following \emph{energy incrementation algorithm} to construct $g$ and $b$.  We shall need two auxiliary $\sigma$-algebras $\B$ and $\B'$, with $\B'$ always being larger than or equal to $\B$.  Also, $\B$ will always be of the form
$\B = \B_{\eps_1, \chi_1} \vee \ldots \vee \B_{\eps_n, \chi_n}$ for some $n$, some $\eps_1,\ldots,\eps_n > 0$, and some $\chi_1,\ldots,\chi_n$, and similarly for $\B'$ (but with different values of $n$); also we will have the bound
\begin{equation}\label{energy-stable}
 \| \E( f|\B')\|_{L^2(\Z/N\Z)}^2 \leq \| \E(f|\B) \|_{L^2(\Z/N\Z)}^2 + \sigma^2/4
 \end{equation}
or equivalently (by Pythagoras' theorem)
\begin{equation}\label{fbb}
 \| \E(f|\B') - \E(f|\B) \|_{L^2(\Z/N\Z)} \leq \sigma/2.
 \end{equation}

\begin{itemize}
\item Step 0: Initialize $\B = \B' = \{0, \Z/N\Z\}$ to be the trivial $\sigma$-algebra.  Note that \eqref{energy-stable} is trivially true at present.
\item Step 1: By construction, we have $\B = \B_{\eps_1, \chi_1} \vee \ldots \vee \B_{\eps_n, \chi_n}$ for some $\eps_1,\ldots,\eps_n > 0$
and linear phase functions $\chi_1,\ldots,\chi_n$.  The function $\E(f|\B)$ is bounded and measurable with respect to $\B$.  By Corollary \ref{sigma-cor} we can thus find $K$ depending on $\delta$, $n$, $\eps_1,\ldots, \eps_n$ such that $\E(f|\B)$ is $(\sigma/2,K)$-almost periodic.
\item Step 2: Set $g := \E(f|\B')$ and $b = f - \E(f|\B')$.  If $\|\hat b\|_{l^\infty} \leq F(\delta,K)$ then we terminate the algorithm; otherwise
we move on to Step 3.
\item Step 3: Since we have not terminated the algorithm, we have $\|\hat b\|_{l^\infty} > F(\delta,K)$.  Using Lemma \ref{sigma-grow}, we
can then find $\eps = F(\delta,K)/C$ and a character $\chi$, with associated $\sigma$-algebra $\B_{\eps,\chi}$, such that
$$ \| \E( b | \B_{\eps,\chi} ) \|_{L^2(\Z/N\Z)} \geq C^{-1} F(\delta,K).$$
From the identity
$$ \E( b | \B_{\eps,\chi} ) = \E( \E( f | \B' \vee \B_{\eps,\chi} ) - \E(f | \B') | \B_{\eps,\chi} )$$
and Pythagoras's theorem, we thus have
$$ \| \E( f | \B' \vee \B_{\eps,\chi} ) - \E(f | \B') \|_{L^2(\Z/N\Z)} \geq C^{-1} F(\delta,K),$$
which by Pythagoras again implies the energy increment
$$ \| \E( f | \B' \vee \B_{\eps,\chi} ) \|_{L^2(\Z/N\Z)}^2 \geq \| \E(f | \B') \|_{L^2(\Z/N\Z)}^2 + C^{-2} F(\delta,K)^2.$$
\item Step 4: We now replace $\B'$ with $\B' \vee \B_{\eps,\chi}$.  If we continue to have the property \eqref{energy-stable},
thne we return to Step 2.  Otherwise, we replace $\B$ with $\B'$ and return to Step 1.
\end{itemize}

Let us first see why this algorithm terminates.  If $\B$ (and hence $K$) is fixed, then each time we pass through Step 4,
the energy $\| \E(f | \B') \|_{L^2(\Z/N\Z)}^2$ increases by at least $C^{-2} F(\delta,K)^2$.  Thus either we terminate the algorithm, or \eqref{energy-stable} must be violated, within $C \sigma^2 / F(\delta,K)^2 = C(F,\delta,K)$ steps.  If the latter occurs, then $B$ is replaced by a new $\sigma$ algebra involving $C(F,\delta,K)$ new characters, with corresponding $\eps$ parameters which are bounded from
below by $C(F,\delta,K)^{-1}$.  This implies that the $K$ quantity associated to $\B$ will be replaced by a quantity of the form
$C(F,\delta,K)$.  Also, the energy $\| \E(f | \B) \|_{L^2(\Z/N\Z)}^2$ will have increased by at least $\sigma^2/4$, thanks to the violation of
\eqref{energy-stable}.  On the other hand, since $f$ was assumed bounded, this energy cannot exceed 1.  Thus we can change $\B$ at most $O(\sigma^{-2})$
times.  Putting all this together we see that ths entire algorithm must terminate in $C(F,\delta)$ steps, and the quantity $K$ will also
not exceed $C(F,\delta)$.  (Note that these constants can be extremely large, as they will involve iterating $F$ repeatedly; however,
the key point is that they do not depend on $N$). 

The claims of the proposition now follow from construction.  Note that $\E(f|\B)$ is $(\sigma/2,K)$-almost periodic by construction,
and hence $g = \E(f|\B')$ will be $(\sigma,K)$-almost periodic thanks to \eqref{fbb}.
\end{proof}

We can now finally prove Roth's theorem.  We let $F: \R^+ \to \R^+ \to \R^+$ be a function to be chosen later, and apply the above
Proposition to decompose $f=g+b$.  By Lemma \ref{ap-recur} we have
$$ \Lambda_3(g,g,g) \geq c(K,\delta) - o_{n,\delta}(1)$$
and then by \eqref{l3-infty} and \eqref{bf} we have
$$ \Lambda_3(f,f,f) \geq c(K,\delta) + O( \delta F(K,\delta) ) - o_{n,\delta}(1).$$
By choosing $F$ sufficiently small, we can absorb the second term in the first, thus
$$ \Lambda_3(f,f,f) \geq c(K,\delta)/2 - o_{n,\delta}(1).$$
Since $K \leq C(F,\delta) = C(\delta)$, the claim \eqref{t3} now follows.
\end{proof}

We remark that there are several other proofs of Roth's theorem in the literature, notably Szemer\'edi's proof based on
density increment arguments and extremely large cubes (see \cite{graham}), and an argument based on the Szemer\'edi regularity lemma (which
in turn requires energy increment arguments in the proof) in \cite{solymosi}.  While these arguments are also important to the theory and both
have generalizations to higher $k$, we will not discuss them here due to lack of space.

\section{Interlude on multilinear operators}

We will shortly turn our attention to Szemer\'edi's theorem.  Based on the preceding section, it is unsurprising that much of the analysis
will revolve around the multilinear form
$$ \Lambda_k(f_0,\ldots,f_{k-1}) := \E( \prod_{j=0}^{k-1} f_j(x + jr) | x, r \in \Z/N\Z )$$
for a large prime $N$.  It turns out that to analyze this multilinear form, it is convenient to generalize substantially and
consider multilinear expressions of the form
\begin{equation}\label{multi-exp}
 \E( K(x) \prod_{j=1}^d F_j(x) | x \in \prod_{j=1}^d A_j )
\end{equation}
where $d \geq 1$ is fixed, $A_1,\ldots,A_d$ are finite non-empty sets,
$K: \prod_{j=1}^d A_j \to \C$ is a fixed kernel, $x = (x_1,\ldots,x_d)$, 
and each $F_i: \prod_{j=1}^d A_j \to \C$ is a bounded function which is independent of the $x_i$ co-ordinate (and thus only depends on the other
$n-1$ co-ordinates).

Henceforth we fix $d$ and $A_1,\ldots,A_d$.
Let $\{0,1\}^d$ be the discrete unit cube.  We need the following notation: if $x^{(0)} = (x^{(0)}_1,\ldots,x^{(0)}_d)$
and $x^{(1)} = (x^{(1)}_1,\ldots,x^{(1)}_d)$ are elements of $\prod_{j=1}^d A_j$, and $\eps = (\eps_1,\ldots,\eps_d) \in \{0,1\}^d$, then we
write $x^{(\eps)} := (x^{(\eps_1)}_1, \ldots, x^{(\eps_d)}_d) \in \prod_{j=1}^d A_j$, and refer to the $2^d$-tuple $( x^{(\eps)} )_{\eps \in \{0,1\}^d}$
of elements in $\prod_{j=1}^d A_j$ as the \emph{cube} generated by $x^{(0)}$ and $x^{(1)}$; this is a cube in the combinatorial sense rather than the geometric sense.  Thus for instance, when $d=2$, the cube generated by $(x,y)$ and $(x',y')$ is
the $4$-tuple consisting of $(x,y)$, $(x,y')$, $(x',y)$, and $(x',y')$.

Now suppose we have a $2^d$-tuple of kernels $K^{(\eps)}: \prod_{j=1}^d A_j \to \C$ for each $\eps \in \{0,1\}^d$.  We define 
the \emph{Gowers inner product} $\langle (K^{(\eps)})_{\eps \in \{0,1\}^d} \rangle_{\Box^d}$ to be 
$$ \langle (K^{(\eps)})_{\eps \in \{0,1\}^d} \rangle_{\Box^d} := \E( \prod_{\eps \in \{0,1\}^d} {\mathcal C}^{|\eps|} K^{(\eps)}(x^{(\eps)}) |
x^{(0)}, x^{(1)} \in \prod_{j=1}^d A_j )$$
where ${\mathcal C} f := \overline{f}$ is the conjugation operator, and $|\eps| := \sum_{j=1}^d \eps_j$.  By separating the $d^{th}$ co-ordinates
of $x^{(0)}$ and $x^{(1)}$, we observe the identity
\begin{equation}\label{ident}
\begin{split} \langle (K^{(\eps)})_{\eps \in \{0,1\}^d} \rangle_{\Box^d}
= \E( 
&\E( \prod_{\underline{\eps} \in \{0,1\}^{d-1}} {\mathcal C}^{|\underline{\eps}|} K^{(\underline{\eps},0)}(\underline{x}^{(\underline{\eps})},y) | y \in A_d ) \\
&{\mathcal C} \E( \prod_{\underline{\eps} \in \{0,1\}^{d-1}} {\mathcal C}^{|\underline{\eps}|} K^{(\underline{\eps},1)}(\underline{x}^{(\underline{\eps})},y) | y \in A_d) | \underline{x}^{(0)}, \underline{x}^{(1)} \in \prod_{j=1}^{d-1} A_j )
\end{split}
\end{equation}
Applying Cauchy-Schwarz in the variables $\underline{x}^{(0)}, \underline{x}^{(1)}$, we conclude that
$$ |\langle (K^{(\eps)})_{\eps \in \{0,1\}^d} \rangle_{\Box^d}|
\leq \langle (K^{(\underline{\eps},0)})_{\eps \in \{0,1\}^d} \rangle_{\Box^d}^{1/2} \langle (K^{(\underline{\eps},1)})_{\eps \in \{0,1\}^d} \rangle_{\Box^d}^{1/2},$$
where $\underline{\eps} := (\eps_1,\ldots,\eps_{d-1})$ are the first $d-1$ co-ordinates of $\eps$; note that \eqref{ident} ensures
that the inner products appearing in the right-hand side of the above equation are non-negative reals.  Of course one has a similar
inequality if we work with the $j^{th}$ co-ordinate instead of the $d^{th}$ co-ordinate for any $1 \leq j \leq d$.  Applying the above
Cauchy-Schwarz inequality once in each co-ordinate, we obtain the \emph{Gowers-Cauchy-Schwarz inequality}
\begin{equation}\label{gcz}
 \langle (K^{(\eps)})_{\eps \in \{0,1\}^d} \rangle_{\Box^d} \leq \prod_{\eps \in \{0,1\}^d} \| K^{(\eps)} \|_{\Box^d}
\end{equation}
where $\| K \|_{\Box^d}$ is the \emph{Gowers cube norm}
$$ \| K \|_{\Box^d} := \langle (K)_{\eps \in \{0,1\}^d} \rangle_{\Box^d}^{1/2^d}.$$
Again, the identity \eqref{ident} ensures that this norm is non-negative.  Using the multilinearity of the Gowers inner product, we then observe for
an arbitrary pair $K_0$, $K_1$ of kernels that
\begin{align*}
\| K_0 + K_1 \|_{\Box^d}^{2^d} &= \langle (K_0 + K_1)_{\eps \in \{0,1\}^d} \rangle_{\Box^d} \\
&= \sum_{A \subset \{0,1\}^d} \langle (K_{1_A(\eps)})_{\eps \in \{0,1\}^d} \rangle_{\Box^d} \\
&\leq \sum_{A \subset \{0,1\}^d} \prod_{\eps \in \{0,1\}^d} \| K_{1_A(\eps)} \|_{\Box^d} \\
&= (\| K_0 \|_{\Box^d} + \| K_1 \|_{\Box^d})^{2^d}
\end{align*}
which thus yields the \emph{Gowers triangle inequality}
$$ \| K_0 + K_1 \|_{\Box^d} \leq \| K_0 \|_{\Box^d} + \| K_1 \|_{\Box^d}.$$
Since the Gowers cube norm is clearly homogeneous, we thus see that $\| \cdot \|_{\Box^d}$ is a semi-norm.  We will later show that it is
in fact a norm when $d \geq 2$; when $d=1$ we have $\| K \|_{\Box^1} = |\E(K(x) | x \in A_1)|$ which is degenerate and thus not a genuine norm.

The significance of the Gowers cube norm to expressions of the form \eqref{multi-exp} lies in the following estimate (which is 
implicit in \cite{ccw} and also in \cite{gowers}).

\begin{lemma}[Van der Corput lemma]\label{vdc-l}  Let $d \geq 1$, let $A_1,\ldots,A_d$ be finite non-empty sets,
let $K: \prod_{j=1}^d A_j \to \C$, and for each $1 \leq i \leq d$ let $F_i: \prod_{j=1}^d A_j \to \C$ is a bounded function which is independent of 
the $x_i$ co-ordinate.  Then we have
$$  |\E( K(x) \prod_{j=1}^d F_j(x) | x \in \prod_{j=1}^d A_j )| \leq \| K \|_{\Box^d}.$$
\end{lemma}

\begin{proof} We induct on $d$.  When $d=1$, the claim becomes
$$  |\E( K(x_1) F_1(x_1) | x_1 \in A_1 )| \leq |\E( K(x_1) | x_1 \in A_1 )|,$$
which follows since $F_1$ is independent of $x_1$ and is bounded.

Now suppose that $d \geq 2$ and the claim has already been proven for $d-1$.  Since $F_d$ is independent of the $x_d$ co-ordinate,
we may abuse notation and interpret $F_d$ as a function on $\prod_{j=1}^{d-1} A_j$ rather than $\prod_{j=1}^d A_d$.  We then
separate off the $x_d$ co-ordinate to write
$$
\E( K(x) \prod_{j=1}^d F_j(x) | x \in \prod_{j=1}^d A_j )
=
\E( F_d(\underline{x}) 
\E( K(\underline{x}, x_d) \prod_{j=1}^{d-1} F_j(\underline{x},x_d) | x_d \in A_d )
| \underline{x} \in \prod_{j=1}^{d-1} A_j) .$$
Since $F_d$ is bounded, we may apply Cauchy-Schwarz in the $\underline{x}$ variable to then obtain
\begin{align*}
|\E( K(x) \prod_{j=1}^d F_j(x) | x \in \prod_{j=1}^d A_j )|
&\leq
\E( 
|\E( K(\underline{x}, x_d) \prod_{j=1}^{d-1} F_j(\underline{x},x_d) | x_d \in A_d )|^2
| \underline{x} \in \prod_{j=1}^{d-1} A_j) \\
&=
\E(
\E( K(\underline{x}, x^{(0)}_d) \overline{K(\underline{x}, x^{(1)}_d)} 
\prod_{j=1}^{d-1} F_j(\underline{x},x^{(0)}_d) \overline{F_j(\underline{x}, x^{(1)}_d)} |
\underline{x} \in \prod_{j=1}^{d-1} A_j) \\
&\quad | x^{(0)}_d, x^{(1)}_d \in A_d )^{1/2}.
\end{align*}
For each fixed $x^{(0)}_d, x^{(1)}_d \in A$ and each $1 \leq j \leq d$, the function
$F_j(\underline{x},x^{(0)}_d) \overline{F_j(\underline{x}, x^{(1)}_d)}$ is a bounded function of
$\underline{x}$.  If we then apply the induction hypothesis we have
$$
|\E( K(x) \prod_{j=1}^d F_j(x) | x \in \prod_{j=1}^d A_j )|
\leq
\E(
\| K(\cdot, x^{(0)}_d) \overline{K(\cdot, x^{(1)}_d)} \|_{\Box^{d-1}} | x^{(0)}_d, x^{(1)}_d \in A_d )^{1/2},$$
so by H\"older's inequality
$$
|\E( K(x) \prod_{j=1}^d F_j(x) | x \in \prod_{j=1}^d A_j )|
\leq
\E(
\| K(\cdot, x^{(0)}_d) \overline{K(\cdot, x^{(1)}_d)} \|_{\Box^{d-1}}^{2^{d-1}} | x^{(0)}_d, x^{(1)}_d \in A_d )^{1/2^d}.$$
But the right-hand side can be re-arranged to be precisely $\|K\|_{\Box^d}$, and the claim follows.
\end{proof}
 
We can now show that $\|\cdot\|_{U^d}$ is a genuine norm when $d \geq 2$:

\begin{corollary} If $d \geq 2$ and $\|K\|_{U^d} = 0$, then $K = 0$.
\end{corollary}

\begin{proof}  Let $(x_1,\ldots,x_d) \in \prod_{j=1}^d A_j$ be arbitrary.  We then define $f_i: \prod_{j=1}^d A_j \to \C$ by defining
$f_i(y_1,\ldots,y_d) = 1$ when $y_j = x_j$ for all $j \neq i$, and $f_i(y_1,\ldots,y_d) = 0$ otherwise.  Applying the previous lemma we thus
see that $K(x_1,\ldots,x_d) = 0$.  Since $(x_1,\ldots,x_d)$ was arbitrary, the claim follows.
\end{proof}

Let us informally call a kernel $K$ \emph{Gowers uniform} if it has small $\Box^d$ norm.  Then the van der Corput lemma then asserts
that Gowers uniform kernels are negligible for the purpose of computing multilinear expressions such as \eqref{multi-exp}.  In particular, when $d=2$,
the $\Box^2$ norm of a kernel $K$ (which can now be interpreted as a linear operator $T_K$ from $L^2(A_1)$ to $L^2(A_2)$) controls the $L^2$ operator norm of $K$.  Indeed, one has the identity 
\begin{equation}\label{ku2}
\begin{split}
\|K\|_{U^2} &= \| T^*_K T_K \|_{HS(L^2(A_1) \to L^2(A_1))}^{1/2}\\
& = 
\| T_K T^*_K \|_{HS(L^2(A_2) \to L^2(A_2))}^{1/2}\\
&= \tr_{A_1}( T^*_K T_K T^*_K T_K )^{1/4} \\
&= \tr_{A_2}( T_K T^*_K T_K T^*_K )^{1/4}
\end{split}
\end{equation}
where $HS$ is the normalized Hilbert-Schmidt norm, and $\tr_A$ is the normalized trace on $A$; equivalently, $\|K\|_{\Box^2}$ is the $l^4$ norm of
the (normalized) singular values of $K$, while the operator norm is the $l^\infty$ norm of these singular values (and the Hilbert-Schmidt norm is
the $l^2$ norm).  Thus one can view the $\Box^d$ norm as a multilinear generalization of the $l^4$ Schatten-von Neumann norm.  This norm has
also arisen in the study of pseudorandom sets and graphs, see for instance \cite{chang}.

Now we specialize to the problem of counting arithmetic progressions in $\Z/N\Z$.

\begin{definition}[Gowers uniformity norm]  Let $f: \Z/N\Z \to \C$ be a function and $d \geq 1$.  Then we define the \emph{Gowers uniformity norm} $\|f\|_{U^d}$
to be the quantity $\| f \|_{U^d} := \| K \|_{\Box^d}$, where $K: (\Z/N\Z)^d \to \C$ is the kernel
$$ K(x_1,\ldots,x_d) := f(x_1+\ldots+x_d).$$
Equivalently, we have
$$ \|f\|_{U^d} := \E( \prod_{\eps \in \{0,1\}^d} {\mathcal C}^{|\eps|} f(x + \sum_{j=1}^d \eps_j h_j)
| x, h_1,\ldots,h_d \in \Z/N\Z )^{1/2^d},$$
or alternatively we have the recursive definitions
\begin{equation}\label{recursive}
 \|f\|_{U^1} := |\E(f)|; \quad \| f \|_{U^{d+1}} := \E( \| f(x+h) \overline{f(x)} \|_{U^d_x}^{2^d} | h \in \Z/N\Z )^{1/2^{d+1}}.
\end{equation}
\end{definition}

Since $\Box^d$ was a norm for $d \geq 2$, we see that $U^d$ is also a norm when $d \geq 2$.  In the $d=2$ case, one can easily verify the
identity
$$ \| f \|_{U^2} := \| \hat f \|_{l^4},$$
which can be viewed as a special case of \eqref{ku2}, observing that the Fourier coefficients of $f$ are essentially the eigenvalues of $K$.
However, for $d \geq 3$ the $U^d$ norm becomes more complicated, and has no particularly useful representation in terms of the Fourier transform.
Using the Gowers-Cauchy-Schwarz inequality, it is possible to show the monotonicity relationship $\| f \|_{U^d} \leq \|f\|_{U^{d+1}}$ for all
$d$; one can also show that $\|f\|_{U^d} \to \|f\|_{L^\infty}$ as $d \to \infty$.  We shall neither prove nor use these facts here.

We can now obtain an analogue of \eqref{lambda4-est}.

\begin{lemma}[Generalized von Neumann theorem]\label{gvn-thm}\cite{gowers}  Let $k \geq 3$, and let $N$ be a prime larger than $k$.  Let
$f_0,\ldots,f_{k-1}$ be bounded functions on $\Z/N\Z$.  Then we have
$$ |\Lambda_k(f_0,\ldots,f_{k-1})| \leq \min_{0 \leq j \leq k-1} \|f_j\|_{U^{k-1}}.$$
\end{lemma}

\begin{proof}  Fix $0 \leq j \leq k-1$; it thus suffices to show that
$$ |\Lambda_k(f_0,\ldots,f_{k-1})| \leq \|f_j\|_{U^{k-1}}.$$
Observe that for any $x_1,\ldots,x_{k-1} \in \Z/N\Z$, the sequence
$$ (x_1 + \ldots + x_{k-1} - (j-i) \sum_{1 \leq i' \leq k: i' \neq j} \frac{1}{j-i'} x_{i'})_{1 \leq i \leq k}$$
is an arithmetic progression of length $k$ in $\Z/N\Z$ (here we are using the hypothesis that $N$ is prime and larger than $k$ in order to
invert $j-i'$).  Conversely, each progression $x, x+r, \ldots, x+(k-1)r$ can be expressed in the above form in exactly the same number of ways
($N^{k-3}$, to be exact).  We may thus write
$$
\Lambda_k(f_0,\ldots,f_{k-1}) = 
\E( \prod_{i=0}^{k-1} f_i(x_1 + \ldots + x_{k-1} - (j-i) \sum_{1 \leq i' \leq k: i' \neq j} \frac{1}{j-i'} x_{i'}) | 
x_1,\ldots,x_{k-1} \in \Z/N\Z ).$$
Now observe that the $i^{th}$ factor in the above sum is bounded and will not depend on $x_i$ when $i \neq j$, and that the $j^{th}$
factor is $f_j(x_1 + \ldots + x_{k-1})$.  Applying the van der Corput lemma and the definition of the $U^{k-1}$ norm, we obtain
the claim.
\end{proof}

Let us informally call a bounded function $f$ \emph{Gowers uniform of order $k-2$} if $\|f\|_{U^{k-1}}$ is small;
thus for instance a function with small $U^2$ norm is linearly uniform, a function with small $U^3$ norm is quadratically uniform, and so forth.
The above lemma then asserts that functions which are Gowers uniform of order $k-2$ have a negligible impact on the $\Lambda_k$ multilinear
form.

\begin{example}\label{u3-ex}  Let $N$ be a prime number, let $P: \Z/N\Z \to \Z/N\Z$ be a polynomial of degree $d$ in the field $\Z/N\Z$, and let
$f(x) := e_N(P(x))$, thus $f$ is a bounded function.  One can easily verify that $\|f\|_{U^{k-1}} = 1$ when $d \leq k-2$
(basically because the $(k-1)^{th}$ derivative of $P$ vanishes), so that $P$ is not uniform of any order $d$ or greater.  (In fact, one
has the more general statement that $\|fg\|_{U^{k-1}} = \|g\|_{U^{k-1}}$ for arbitrary $g$ and whenever $d \leq k-2$; thus the $U^{k-1}$ norm is invariant under polynomial phase modulations of degree $k-2$ or less).  On the other hand, 
one can verify that $\| f \|_{U^{k-1}} = O_d(N^{-1/2})$ when $d > k-2$; this is easiest to accomplish when $d = k-1$, and the remaining cases follow
by monotonicity (or van der Corput type arguments for Weyl sums).  Thus $P$ is uniform of order $d-1$ or less.  The intuition to have here is that
a bounded function is (heuristically) uniform of order $d$ iff its phase is ``orthogonal'' to all polynmial phases of degree $d$ or less. 
In the $d=1$ case this intuition is precise: linear uniformity corresponds to being orthogonal to linear phase functions, as
the estimates \eqref{lest} already attest to.  When $d \geq 2$ however this intuition is harder to pin down, and the theory
is still not completely understood.  

Now consider a quadratic polynomial $P(x)$, with corresponding quadratic phase function $f(x) := e_N(P(x))$.  From the identity
$$ P(x) - 3P(x+r) + 3P(x+2r) - P(x+3r) = 0$$
(which reflects the fact that the third derivative of $P$), we observe that
$$ \Lambda_3(f, \overline{f}, f, \overline{f}) = 1.$$
Thus $f$ is non-negligible for the purposes of computing the $\Lambda_3$ form.  This is despite $f$ being linearly uniform (all the Fourier
coefficients of $f$ is $O(N^{-1/2})$, as one sees from the classical theory of Gauss sums).  This shows that for the purposes of analyzing
$\Lambda_3$, it is really quadratic uniformity which is the concept to be studied, not linear uniformity.  Similarly, the concept of
being Gowers uniform of order $k-2$ is the one which is related to the form $\Lambda_k$, which in turn counts arithmetic progressions of
length $k$.
\end{example}

\section{Progressions of length 4}

With the above machinery, we can now sketch two different proofs of Szem\'eredi's theorem for progressions
of length 4.  (These arguments also extend, with some additional difficulties, to higher $k$, but we will not discuss
these technicalities here).  The first proof we present is due to Gowers \cite{gowers-4} and can be viewed as a generalization
of Roth's Fourier-analytic argument, being a density-incrementation argument using quadratic Fourier analysis instead of 
linear Fourier analysis.  The second proof is adapted from that in \cite{tao:szemeredi}, which in turn is
based on the original ergodic theory arguments of Furstenberg and co-authors \cite{furst}, \cite{FKO}.  
It is a generalization of the second proof of Roth's theorem given earlier; in particular, it is
is an energy-incrementation argument based on the decomposition of an arbitrary function into a ``almost periodic
function of order 2'' and a quadratically uniform function.

We begin by discussing Gowers' proof, though we shall omit many of the details which pertain to arithmetic combinatorics.  
Once again, we have a subset $A$ of $[1,N]$, which we embed into a cyclic group $\Z/p\Z$ of prime order.
We split $f = 1_A = g + b$, where $g = \E(f)$ and $b = f - \E(f)$.  If $b$ is quadratically uniform in the sense that $\|b\|_{U^3}$ is
suitably small (less than $c \delta^C$ for some absolute constants $c, C > 0$) then, by using Lemma \ref{gvn-thm} to develop an analogue of
Proposition \ref{gbgb}, then one can easily obtain non-trivial lower bounds for $\Lambda_4(f,f,f,f)$ and thus establish plenty of
arithmetic progressions of length 4 in $A$.

The difficulty comes in the ``hard case'', when $b$ is not quadratically uniform, so that $\|b\|_{U^3}$ is relatively large.  The difficulty
here is that unlike the $U^2$ norm, which is the $l^4$ norm of the Fourier transform, the $U^3$ norm is not easily related to the Fourier transform;
for instance in Example \ref{u3-ex} we saw that there were functions which had very small Fourier transform but had large $U^3$ norm.  Nevertheless, it is still possible to use this information to deduce some structural information about $A$. The
situation can be clarified somewhat by considering a model problem, which is to determine all functions of the form $b = e_p(\phi(x))$
which had the maximal $U^3$ norm of $1$, where $\phi: \Z/p\Z \to \Z/p\Z$ is a phase function.  Expanding out the $U^3$ norm, we see that this is equivalent to asking that
\begin{equation}\label{phi3}
\phi(x+r+s+t) - \phi(x+r+s) - \phi(x+r+t) - \phi(x+s+t) + \phi(x+r) + \phi(x+s) + \phi(x+t) - \phi(x) = 0
\end{equation}
for all $x,r,s,t \in \Z/p\Z$.  This is an ``arithmetic'' way of asserting that the third derivative of $\phi$ vanishes.  It in fact implies that $\phi$ is a quadratic polynomial, $\phi(x) = a x^2 + b x + c$ (whereas in contrast, the assertion that $b$ would have a maximal Fourier coefficient of 1 is equivalent to asserting that $\phi$ is a linear polynomial).  To see this, let us adopt the notation that for any function $f: \Z/p\Z \to \R/\Z$ and any shift $h \in \Z/p\Z$, that $f_h: \Z/p\Z \to \R/\Z$ denotes the ``derivative'' $f_h(x) := f(x+h) - f(x)$.  Then we have
\begin{equation}\label{phih}
\phi_h(x+s+t) - \phi_h(x+s) - \phi_h(x+t) + \phi_h(x) = 0 \hbox{ for all } x,s,t \in \Z/p\Z.
\end{equation}
It is easy to see that this implies that $\phi_h$ is linear, i.e. we have 
\begin{equation}\label{phi-ab}
\phi_h(x) = a(h) x + c(h) \hbox{ for all } x,h \in \Z/p\Z
\end{equation}
for some $a(h)$, $c(h) \in \Z/p\Z$.
(This is easiest seen by first subtracting $\phi_h(0)$ from $\phi_h$, at which point $\phi_h$ becomes additive).  
To conclude from this that $\phi$ is quadratic, one would need to firstly show that $a(h)$ and $c(h)$ have some linearity properties in $h$, and then ``integrate'' the
equation \eqref{phi-ab} to obtain a quadratic expression for $\phi(x)$.

To attain these goals, we rewrite \eqref{phi-ab} as the functional equation
\begin{equation}\label{ahunk}
a(h) x = \phi(x+h) - \phi(x) - c(h) \hbox{ for all } x,h \in \Z/p\Z.
\end{equation}
We can isolate $a(h)$ in this equation by taking suitable ``derivatives''.  For instance, if one replaces $x$ by $x+s$ in the above formula to obtain
$$ a(h) (x+s) = \phi(x+s+h) - \phi(x+s) - c(h) \hbox{ for all } x,h,s \in \Z/p\Z.$$
and then subtracts the two equations, one obtains
\begin{equation}\label{ahs}
 a(h) s = \phi_s(x+h) - \phi_s(x) \hbox{ for all } x,h,s \in \Z/p\Z
\end{equation}
thus eliminating the unknown function $c(h)$.  Similarly, by replacing $h$ by $h+t$ and then subtracting, we can eliminate the $\phi_s(x)$ term
to obtain
$$ a_t(h) s = \phi_{st}(x+h) \hbox{ for all } x,h,s,t \in \Z/p\Z.$$
Finally, by replacing $x,h$ by $x-u,h+u$ and subtracting again to
eliminate the $\phi_{st}(x+h)$ term, one obtains
\begin{equation}\label{atuh}
 a_{tu}(h) s = 0 \hbox{ for all } h,s,t,u \in \Z/p\Z
\end{equation}
and thus $a$ obeys the functional equation
\begin{equation}\label{atuh-linear}
a(h+t+u) - a(h+t) - a(h+u) + a(h) = 0 \hbox{ for all } h,u,t \in \Z/p\Z
\end{equation}
which as observed earlier implies that $a(h)$ is linear, thus 
\begin{equation}\label{a-alphabeta}
a(h) = \alpha h + \beta \hbox{ for some } \alpha,\beta \in \Z/p\Z.
\end{equation}
(One can in fact force $\beta$ to equal zero, basically because $a(0) = 0$, but we will not do so here).
Now the function $\alpha h x$ can be explicitly integrated (modulo a lower order term) using the quadratic primitive 
\begin{equation}\label{primitive-def}
F(x) := \frac{\alpha}{2} x^2,
\end{equation}
in the sense that $F_h(x) = \alpha h x + \frac{\alpha}{2} h^2$.  Thus if we define $\phi'(x) := \phi(x) - F(x)$ and $\phi''(x) := \phi(x) - F(x) + \beta x$, then by \eqref{phi-ab}, $\tilde \phi$
obeys the functional equation
\begin{equation}\label{phip}
\phi'(x+h) - \phi''(x) = \beta x + c(h) - \frac{\alpha}{2} h^2 \hbox{ for all } x, h \in\Z/p\Z.
\end{equation}
Replacing $x$ by $x+k$ and subtracting, we obtain that
$$ \phi'(x+h+k) - \phi'(x+h) - \phi''(x+k) + \phi''(x) = 0 \hbox{ for all } x,h,k \in \Z/p\Z$$
which then implies that $\phi'$ and $\phi''$ is linear.  Since $\phi = \phi' + F$, we thus see that $\phi$ is quadratic as claimed.

This concludes the treatment of the model problem.  Thanks to the work of Gowers \cite{gowers-4}, it turns out that the general strategy used to solve this model problem can also be used to handle the general case.  Indeed, if a function $b$ has large $U^3$ norm (where by ``large'' we mean
``larger than $C^{-1} \delta^C$ for some absolute constant $C > 0$''), then by \eqref{recursive} the function
$b(x+h) \overline{b(x)}$ will have large $U^2$ norm for a large percentage of $h \in \Z/p\Z$ (this is the analogue of \eqref{phih}).
Since $U^2$ norms imply large Fourier coefficients, we thus see that for all $h$ in a large fraction $H \subset \Z/p\Z$ of $\Z/p\Z$ we can find
$a(h), c(h) \in \Z/p\Z$ such that
\begin{equation}\label{buh}
\Re \E( b(x+h) \overline{b(x)} e_p( a(h) x + c(h) ) | x \in \Z/p\Z ) \geq C^{-1} \delta^C
\end{equation}
and hence
$$ |\E( b(x+h) \overline{b(x)} e_p( -(a(h) x + c(h)) ) 1_H(h) | x,h \in \Z/p\Z )| \geq C^{-1} \delta^C.$$
As with the model problem, the task would now be to obtain some linearity control on $a$.  This can be obtained by a Cauchy-Schwarz
argument; there are a number of permutations of this argument, but we shall give one which is based on the van der Corput lemma, 
Lemma \ref{vdc-l}.  Let us first change variables
$x = y_1 - y_2$, $h = y_2 - y_3$ to obtain
$$ |\E( b(y_1 - y_3) \overline{b(y_1-y_2)} 1_H(y_2-y_3) e_p(-c(y_2-y_3)) K(y_1,y_2,y_3) | y_1,y_2,y_3 \in \Z/p\Z )| \geq C^{-1} \delta^C,$$
where 
$$ K(y_1,y_2,y_3) := e_p( - a(y_2-y_3) (y_1-y_2) ) 1_H(y_2-y_3).$$
If we then apply Lemma \ref{vdc-l}, we conclude that
$$ \| K \|_{\Box^3} \geq C^{-1} \delta^C.$$
Raising this to the eighth power and expanding out the left-hand side, one eventually obtains (after some change of variables)
$$ \E( e_p( -(a(h+t+u) - a(h+t) - a(h+u) + a(h)) s ) 
1_{h,h+u,h+t,h+t+u \in H}| h,t,s,u \in \Z/p\Z ) \geq C^{-1} \delta^C$$
(this is the analogue of \eqref{atuh}).  The average in $s$ can be computed explicitly, and we then obtain
\begin{equation}\label{atuh-avg}
\P( h,h+t,h+u,h+t+u \in H; a(h+t+u) - a(h+t) - a(h+u) + a(h) = 0 | h,t,u \in \Z/p\Z ) \geq C^{-1} \delta^C
\end{equation}
(cf. \eqref{atuh-linear}).  This is now a purely arithmetic-combinatorial statement about $a$, involving no oscillation; it says
that $a$ behaves like an (affine-)linear function ``a significant fraction of the time''.  In analogy with \eqref{a-alphabeta}
It is then tempting to conjecture from this that $a(h)$ should in fact \emph{equal} an affine linear function $\alpha h + \beta$ for a significant fraction of the time, i.e. we should be able to find $\alpha, \beta \in \Z/p\Z$ such that
\begin{equation}\label{aab-avg}
\P( h \in H; a(h) = \alpha h + \beta | h \in \Z/p\Z ) \geq C^{-1} \delta^C
\end{equation}
(note that in the converse direction, that one can use \eqref{aab-avg} and a Cauchy-Schwarz argument to obtain \eqref{atuh-avg}).
Suppose for the moment that one could indeed deduce \eqref{aab-avg} from \eqref{atuh-avg}.  Then we can introduce the primitive function
\eqref{primitive-def} as before, and define $b'(x) := b(x) e_p( -F(x))$ and
$b''(x) := b'(x) e_p(\beta x)$; we then see from \eqref{buh} that
$$
\Re \E( b'(x+h) \overline{b''(x)} e_p( \frac{\alpha}{2} h^2 - c(h) ) | x \in \Z/p\Z ) \geq C^{-1} \delta^C
$$
for all $h \in H$ (cf. \eqref{phip}).  In particular we see that
$$ |\E( b'(x+h) \overline{b''(x)}  | x \in \Z/p\Z )| \geq C^{-1} \delta^C.$$
Taking $L^2$ norms of both sides and using Plancherel, we obtain
$$ \| \widehat{b'} \widehat{b''} \|_{l^2} \geq C^{-1} \delta^C,$$
and thus by H\"older's inequality
$$ \| b' \|_{U^2}^4  \geq C^{-1} \delta^C.$$
To summarize, we started with a function $b$ with large $U^3$ norm, and then were able to locate a quadratic modulation $b'$ of $b$ which
in fact had large $U^2$ norm.  Since we already know that a large $U^2$ norm would imply a large Fourier coefficient, we could thus deduce the existence
of a $\xi \in \Z/p\Z$ such that $\widehat{b'}(\xi)$ is large, which would then imply that the original function $b$ had large correlation
with a quadratic phase function $\chi(x) := e_p( P(x) )$ for some quadratic polynomial $P: \Z/p\Z \to \Z/p\Z$, thus
$|\langle b, \chi \rangle| \geq C^{-1} \delta^C$.  One can now proceed as in the density increment proof of Roth's theorem, but with the Bohr
sets in $\B_{\eps,\chi}$ now being replaced by ``quadratic Bohr sets'' .  This eventually gives us a density increment of the form
\eqref{pnaq} on a quadratic Bohr set $\chi^{-1}(Q)$; one can then use Weyl's theorem on equidistribution of quadratic polynomials mod $p$ to locate a
reasonably long arithmetic progression (of length at least $c N^c$ for some absolute constant $c > 0$, if $N$ is sufficiently large
depending on $\delta$) on which one has a density increment, at which point we may repeat Roth's argument.  We omit the details, referring the
reader instead to \cite{gowers-4}.  

We return briefly now to a step glossed over in the above sketch, namely the deduction of \eqref{aab-avg} from \eqref{atuh-avg}.  As it turns
out, this implication is false as stated; it is possible for $a$ to be additive in the sense of \eqref{atuh-avg} without being approximately
linear in the sense of \eqref{aab-avg}, because $a$ may instead be behaving like a ``higher-dimensional'' linear function.  An example of 
this is as follows. Let $M$ be an integer between $\sqrt{p}/4$ and $\sqrt{p}/2$, let $H := \{ n + 2Mm: 1 \leq n,m \leq M \}$, and let $a: H \to \Z/p\Z$ be the function $a(n + 2Mm) = \alpha n + \beta m$ for some fixed $\alpha, \beta \in \Z/p\Z$.  Then one can easily verify that $a$ obeys the property \eqref{atuh-avg} but not \eqref{aab-avg} (if $\beta \neq 2M\alpha$).  The set $H$ is an example of a \emph{two-dimensional arithmetic progression},
and the function $a$ given here is a \emph{generalized linear function} on this progression; more generally one can define the notion of
a generalized arithmetic progression (of arbitrary dimension), and of a generalized linear function on this progression; it is possible then to
obtain a deduction of the form \eqref{atuh-avg} $\implies$ \eqref{aab-avg} but with the role of $\alpha h + \beta$ being played by these
generalized linear functions; also, for technical reasons (having to do with relatively poor constants in a certain inverse theorem from additive combinatorics known as Freiman's theorem) one must with the lower bound of $C^{-1} \delta^C$ by a smaller quantity such as
$\exp( - C \delta^{-C})$; it is not known whether this exponential loss has to be removed.  The deduction here requires a combination of techniques
from combinatorial graph theory, probabilistic combinatorics, Fourier analysis, and the geometry of lattices and Bohr sets; it is somewhat involved and
we will not go into the details here, referring the reader instead to \cite{gowers-4}.  

The remainder of Gowers' argument in \cite{gowers-4} is concerned with how to use the fact that $a$ is approximately equal to a higher-dimensional
linear function to again deduce a density increment of $A$ on some sub-progression.  This is again done mainly by Weyl's theory of uniform distribution; however in \cite{gt-qm} an alternate argument was developed, which is based on locating a primitive $F$ to $a$.  This argument closely mimics the one given in the one-dimensional case when $a(h) \approx \alpha h + \beta$; however, there is an additional difficulty in the
higher-dimensional case, namely that not every linear function has a primitive; instead, only the ``self-adjoint'' linear functions do.  
This has to do with the fact that quadratic forms in higher dimensions (the analogue of quadratic polynomials in one dimension) are associated to symmetric matrices rather than general matrices.  Fortunately, one can show that the function $a$ does indeed obey the required symmetry property.
Rather than give the precise statement and proof of this assertion in detail, we sketch how it works in a model case.  Here we consider solutions to the equation \eqref{ahunk},
but now $x, h$ take values in a vector space $V := (\Z/p\Z)^n$, and $a(h)$ is now a linear transformation from $V$ to $\Z/p\Z$.  By arguing as before,
we conclude that $a(h) = \alpha h + \beta$, where $\beta$ is now a linear transformation from $V$ to $\Z/p\Z$, and $\alpha$ is a bilinear
form from $V \times V$ to $\Z/p\Z$.  Inserting this back into \eqref{ahs}, we obtain
$$
\alpha(h,s) + \beta s = \phi(x+h+s) - \phi(x+h) - \phi(x+s) + \phi(x) \hbox{ for all } x,h,s \in V.
$$
Now we proceed a little differently to before.  If we replace $h, s$ by $h+u, s-u$ and subtract, we obtain
$$
\alpha(h+u,s-u) - \alpha(h,s) - \beta u = - \phi_u(x+h) - \phi_{-u}(x+s) \hbox{ for all } x,h,s,u \in V.
$$
If now we replace $x,h,s$ by $x+t,h-t,s-t$ and subtract, we obtain
$$
\alpha(h+u-t,s-u-t) - \alpha(h-t,u-t) - \alpha(h+u,s-u) + \alpha(h,s) = 0 \hbox{ for all } h,s,u,t \in V.
$$
Using the bilinearity of $\alpha$, this simplifies to
$$ \alpha(t,u) - \alpha(u,t) = 0 \hbox{ for all } u,t \in V$$
which shows that $\alpha$ is symmetric.  In particular this allows us to construct a primitive $F$ by the formula
$F(x) := \frac{1}{2} \alpha(x,x)$, and the previous argument now proceeds as before.  Back in the original setting of
a function $b$ with large $U^3$ norm, an analogous argument allows us to locate a ``generalized quadratic polynomial phase function''
$\chi(x) := e_p( P(x) )$ such that $\langle b, \chi \rangle$ is somewhat large; see \cite{gt-qm} for a rigorous statement and proof
of this ``inverse theorem for the $U^3$ norm''.  (Interestingly, there are some closely related results arising from ergodic theory; see \cite{host-kra2}, \cite{ziegler}).

This concludes our discussion of Gowers' proof of Szemer\'edi's theorem for progressions of length 4; the argument also extends to higher $k$
(see \cite{gowers}) though with some non-trivial additional difficulties; also, it is not at present clear whether the higher $U^d$ norms also
enjoy an inverse theorem.  We now briefly discuss another proof of this theorem, which extends the energy increment proof for progressions of length three discussed earlier.  There are many proofs in this spirit, starting with the work of Furstenberg \cite{furst}, \cite{FKO} (and a related
energy-incrementation argument also appears in \cite{szemeredi}); we shall loosely follow the version of this argument from \cite{tao:szemeredi}.
For sake of simplicity we shall confine our discussion to the $k=4$ case only.

As it turns out, large portions of the energy increment proof generalize without difficulty to obtain progressions of arbitrary length.  The main
difficulty is to replace the concept of an $(\delta,K)$-almost periodic function with a ``higher order'' generalization.  The definition given in
Definition \ref{ap-def} relies too heavily on linear phase functions, and we have already seen some difficulties in extending that concept
to higher orders; for instance, we still do not have a satisfactory theory of what a ``quadratically quasiperiodic function''
should be, although there are some very promising developments in the ergodic theory of nilfactors (see e.g. \cite{host-kra2}, \cite{ziegler}, \cite{ziegler2}) which should shed light on this question very soon.  However, it is well understood by now how to generalize the more
general concept of an almost periodic function.  In ergodic theory, a function $f$ in a measure-preserving system $(X, {\mathcal B}, \mu, T)$ is
said to be almost periodic if the orbit $\{ T^n f : n \in \Z \}$ is precompact, and in particular can be approximated to arbitrary accuracy by a subset of a finite-dimensional space.  In the discrete setting of $\Z/N\Z$, every function is periodic of order $N$ and is thus, technically speaking, every function is almost periodic.  However one can still extract a useful concept of almost periodicity by making the concept of ``precompact'' more quantitative.  One such way of doing so is

\begin{definition}[Uniform almost periodicity norms]\label{uap-def}\cite{tao:szemeredi}  If $A$ is a shift-invariant Banach algebra of functions on $\Z_N$, we define
the space $UAP[A]$ to be the space of all functions 
$F$ for which the orbit $\{ T^n F: n \in \Z \}$ has a representation of the form
\begin{equation}\label{representation}
 T^n F = M \E( c_{n,h} g_h ) \hbox{ for all } n \in \Z_N
\end{equation}
where $M \geq 0$, $H$ is a finite non-empty set, $g = (g_h)_{h \in H}$ is a collection of bounded functions, $c = (c_{n,h})_{n \in \Z_N, h \in H}$ is a collection of functions in $A$ with $\|c_{n,h}\|_{A} \leq 1$, 
and $h$ is a random variable taking values 
in $H$.  We define the norm $\|F\|_{UAP[A]}$ to be the infimum of $M$ over all possible representations
of this form.  
\end{definition}

The formula \eqref{representation} is a quantitative assertion that the orbit $\{ T^n F: n \in \Z\}$ can be represented efficiently by what is essentially a finite-dimensional approxmation, and is thus an assertion of precompactness ``relative to $A$''.
It can be shown (see \cite{tao:szemeredi}) that $UAP[A]$ is a shift-invariant Banach algebra.
If we let $A$ be the trivial Banach algebra of constant functions (so that the $c_{n,h}$ are constants, with $\|c_{n,h}\|_A = |c_{n,h}|$) then we abbreviate $UAP[A]$ as $UAP^1$, and refer to functions with bounded $UAP^1$ norm as \emph{linearly uniformly almost periodic}.  For instance, one can show
that any $K$-quasiperiodic function is linearly uniformly almost periodic, with a $UAP^1$ norm of at most $K$.  In particular, linear phase functions
are linearly uniformly almost periodic, with a $UAP^1$ norm of exactly 1.

One can then define the space $UAP^2 := UAP[UAP^1]$ of quadratically uniformly functions, which are roughly speaking the space of functions which
are almost periodic relative to the linearly almost periodic functions.  For example, consider the function $f(x) := e_N(x^2)$.  This function is very far from being linearly almost periodic - in the sense that the $UAP^1$ norm is huge - because the translates $T^n f(x) = e_N(x^2 + 2nx + n^2)$ are
all quite distinct and cannot efficiently be expressed as linear combinations of a small number of functions.  On the other hand, we may write
$T^n f = c_n g$ where $g := f$ and $c_n(x) = e_N(2nx + n^2)$, and note that each $c_n$, being a linear phase function, lies in $UAP^1$ with small norm.
Thus this function is quadratically almost periodic; in fact, it lies in $UAP^2$ with norm 1.  The property of being quadratically almost periodic
strictly generalizes the concept of a \emph{quadratic eigenfunction} in ergodic theory; see e.g. \cite{ziegler}, \cite{ziegler2} for further
discussion.

The concept of quadratic almost periodicity (bounded $UAP^2$ norm) is in many ways dual to that of quadratic uniformity (small $U^3$ norm).  We present
three results supporting this claim.  The first is the duality inequality
$$ |\langle f, F \rangle| \leq \| f \|_{U^3} \|F\|_{UAP^2},$$
which can be proven by a simple Cauchy-Schwarz argument, see \cite{tao:szemeredi}.  Secondly, if $f$ is such that $\|f\|_{U^3}, \|f\|_{L^\infty} \leq 1$, and we let ${\mathcal D} f$ denote
the \emph{dual function}
$$ {\mathcal D} f(x) := \E( \overline{f(x+a) f(x+b) f(x+c)} f(x+a+b) f(x+a+c) f(x+b+c) \overline{f(x+a+b+c)} | a,b,c \in \Z_N )$$
then ${\mathcal D} f$ lies in $UAP^2$ with a norm of at most 1; again, see \cite{tao:szemeredi}.  Furthermore, we have the correlation identity
$$ \langle f, {\mathcal D} f \rangle = \| f \|_{U^3}^8.$$
By using these dual function to replace the role of linear (or quadratic) phase functions, one can obtain the following 
variant of Proposition \ref{qkn}:

\begin{proposition}[Quantitative Koopman-von Neumann theorem]\cite{tao:szemeredi}  Let $F: \R^+ \times \R^+ \times \R^+ \to \R^+$ be an arbitrary function, let $0 < \sigma < \delta \leq 1$, and let $f$ be any bounded non-negative
function on $\Z/N\Z$ obeying \eqref{mean}.  Then there exists a quantity $0 < K \leq C(F, \delta, \sigma)$ and a decomposition $f = g+b$, where $g$ is bounded, non-negative, has mean $\E(g)= \E(f)$, and we have the bound
$$
\| b \|_{U^3} \leq F(\delta,\sigma,K).
$$
Furthermore we have an additional decomposition $g = \tilde g + e$ with $\tilde g$ non-negative and the bounds
$$ \| \tilde g \|_{UAP^2} \leq K; \quad \| e \|_{L^2} \leq \sigma.$$
\end{proposition}  

The proof of this Proposition proceeds by an energy incrementation argument very similar to Proposition \ref{qkn}; one begins with the trivial splitting
$f = \E(f) + (f-\E(f))$, and whenever the bad function $b$ fails to be quadratically uniform, one uses the dual function ${\mathcal D} b$ (which is 
quadratically almost periodic) to refine the $\sigma$-algebra used to construct the good function $g$, thus increasing the energy of $g$ by a non-trivial amount.

By combining this with the generalized von Neumann theorem in Lemma \ref{gvn-thm}, we can conclude the proof of Szemer\'edi's theorem in
this case once we show the analogue of Lemma \ref{ap-recur}:

\begin{theorem}[Almost periodic functions are recurrent]\label{recurrence}\label{tao:szemeredi}  Let $g, \tilde g$ be non-negative bounded functions such that we have the estimates
\begin{align}
 \| g - \tilde g \|_{L^2} &\leq \frac{\delta^2}{4096}\label{central-0}\\
 \E( g | \Z_N ) &\geq \delta \label{central-1}\\
 \| \tilde g \|_{UAP^2} &< M\label{central-2}
\end{align}
for some $0 < \delta, M < \infty$.
Then we have
\begin{equation}\label{cdm}
 \E( g(x) T^n g(x) T^{2n} g(x) T^{3n} g(x) | x,r \in \Z_N) \geq c(\delta,M)
\end{equation}
for some $c_0(\delta,M) > 0$.
\end{theorem}

The proof of this theorem is the most difficult component of the argument; it uses the uniform almost periodicity control on $\tilde g$
to ``color'' the orbit of $T^n \tilde g$ and hence $T^n g$, and then invokes the van der Waerden theorem \cite{vdw} to extract arithmetic progressions
from $g$.  As such, this part of the argument can be considered to be more combinatorial than ergodic or analytic in nature.

\section{Progressions in the primes}

There are many questions concerning the distribution of the prime numbers (and of various configurations of 
prime numbers), which has motivated a large portion of analytic number theory.  One of the basic results
in the subject is of course the \emph{prime number theorem}, which asserts that the number of primes between 1 and $N$
asymptotically approaches $N/\log N$ as $N \to \infty$, or in other words
$$ \# \{ 1 \leq n \leq N: n \hbox{ is prime} \} = \frac{N}{\log N}(1 + o(1)),$$
where we use $o(1)$ to denote a quantity which goes to zero as $N \to \infty$.  

It is convenient to normalize the prime number theorem in a different form.  Define the \emph{von Mangoldt function} $\Lambda: \Z^+ \to \R$
by setting $\Lambda(n) := \log p$ whenever $n = p^j$ is a power of a prime $p$ for some $j \geq 1$, and $\Lambda(n) = 0$ otherwise; the significance
of this function to number theory lies in the identity
\begin{equation}\label{logn}
 \log n = \sum_{d | n} \Lambda(d)
 \end{equation}
for all integers $n$ (where the sum is over all integers $d$ dividing $n$), which is a restatement of the unique factorization theorem.
The Von Mangoldt function is essentially supported on the primes (there are also the squares and higher powers of primes, but they are 
extremely sparse, and in practice are completely negligible, contributing only to the $o(1)$ error terms). Then the prime 
number theorem is easily seen to be equivalent to
$$ \frac{1}{N} \sum_{1 \leq n \leq N} \Lambda(n) = 1 + o(1).$$
The expression on the left-hand side can be viewed as an average or \emph{expectation} for $\Lambda$; we shall emphasize this probabilistic
(or ergodic) perspective by writing it as $\E( \Lambda(n): 1 \leq n \leq N)$; more generally, we write 
$\E( f(n) | n \in A)$ for $\frac{1}{|A|} \sum_{n \in A} f(n)$
whenever $A$ is a finite set.  Thus $\Lambda$ has an average value of $1 + o(1)$.  The error can be improved; for instance the famous \emph{Riemann
hypothesis} is equivalent to the claim
$$ \E( \Lambda(n) | 1 \leq n \leq N ) = 1 + O( N^{-1/2} \log^2 N).$$
However the improved error estimates are not central to the results we shall discuss here, which are in some sense more focused on the main term
in such estimates involving the primes.

Now we consider how to count other patterns inside the primes.  One of the oldest (and still unsolved) problems in the field is the \emph{twin prime
conjecture}, which asks whether there are an infinite number of primes $p$ such that $p+2$ is also prime.  This would be implied by
the statement
$$ \liminf_{N \to \infty} \E( \Lambda(n) \Lambda(n+2): 1 \leq n \leq N ) > 0.$$
is non-zero for infinitely many $N$.  In fact Hardy and Littlewood made the stronger conjecture, the \emph{Hardy-Littlewood prime tuple conjecture}
\cite{hardy-littlewood},
which would imply the twin prime conjecture, and would indeed verify the stronger estimate
$$ \E( \Lambda_N(n) \Lambda_N(n+2): 1 \leq n \leq N ) = B_2 + o(1)$$
where $B_2$ is the \emph{Twin prime constant} 
\begin{align*}
B_2 &:= \prod_p \frac{\P( n, n+2 \hbox{ coprime to } p | n \in \Z/p\Z )}{\P( n \hbox{ coprime to } p | n \in \Z/p\Z) \P( n+2 \hbox{ coprime to } p | n \in \Z/p\Z)} \\
&= 2 \prod_{p > 3} \frac{ p(p-2) }{(p-1)^2} \\
&= 1.32032\ldots
\end{align*}
A related problem is the \emph{strong Goldbach conjecture} - whether every even number (larger than 4) can
be written as the sum of two primes; this is essentially the same as asking whether
$$ \E( \Lambda(n_1) \Lambda(n_2): 1 \leq n_1, n_2 \leq N; n_1+n_2 = N )$$
is non-zero for all even integers $N$.  The Hardy-Littlewood prime tuple conjecture here would imply
that
$$ \E( \Lambda(n_1) \Lambda(n_2): 1 \leq n_1, n_2 \leq N; n_1+n_2 = N ) = G_2(N) + o(1)$$
where 
\begin{align*}
G_2(N) &:= \prod_p \frac{\P( n_1,n_2 \hbox{ coprime to } p | n_1,n_2 \in \Z/p\Z; n_1+n_2=N )}
{\prod_{j=1}^2 \P( n_j \hbox{ coprime to } p | n_1,n_2 \in \Z/p\Z; n_1+n_2 = N)} 
\end{align*}
which vanishes when $N$ is odd, and is equal to
$$ G_2(N) = B_2 \prod_{p|N; p \geq 3} \frac{p-1}{p-2} \geq B_2 > 0$$
when $N$ is even.  Thus the prime tuple conjecture would imply the strong Goldbach conjecture for sufficiently large $N$.

The \emph{weak Goldbach conjecture}, which is essentially proven (thanks primarily to the work of Vinogradov \cite{vinogradov}), asserts that every odd number $N$ larger than $5$ can be written as the sum of three primes.  (By ``essentially proven'' I mean that this conjecture has been verified for $N \leq 10^{17}$ and also rigourously proven for $N \geq 10^{43000}$).  This is essentially asking for the quantity
$$ \E( \Lambda(n_1) \Lambda(n_2) \Lambda(n_3): 1 \leq n_1,n_2,n_3 \leq N; n_1+n_2+n_3 = N )$$
to be positive for all odd integers $N$.  The work of Vinogradov implies
$$ \E( \Lambda(n_1) \Lambda(n_2) \Lambda(n_3): 1 \leq n_1,n_2,n_3 \leq N; n_1+n_2+n_3 = N ) = G_3(N) + o(1)$$
where 
\begin{align*}
G_3(N) &:= \prod_p \frac{\P( n_1,n_2,n_3 \hbox{ coprime to } p | n_1,n_2,n_3 \in \Z/p\Z; n_1+n_2+n_3=N )}
{\prod_{j=1}^3 \P( n_j \hbox{ coprime to } p | n_1,n_2,n_3 \in \Z/p\Z; n_1+n_2+n_3 = N)}.
\end{align*}
This quantity is positive and bounded away from zero for all odd $N$; thus Vinogradov's work implies the weak Goldbach conjecture for all
sufficiently large $N$; to resolve the remaining cases it is thus natural to try to sharpen the $o(1)$ error term.  (For instance, the weak Goldbach
conjecture is known to be true if one assumes the generalized Riemann hypothesis, which is extremely useful in improving these error terms).
One can generalize Vinogradov's result to sums of $k$ primes for any $k \geq 3$; but as we shall explain later, the $k=2$ case is much more difficult
and well beyond the reach of existing techniques.

Now we turn to arithmetic progressions in the primes.  In 1933 van der Corput \cite{van-der-corput} (see also \cite{chowla}) established that the primes contain infinitely many arithmetic progressions of length 3; indeed we know the
significantly stronger statement that the Hardy-Littlewood conjecture holds in this case, or more explicitly that
\begin{equation}\label{vdc}
\E( \Lambda(n) \Lambda(n+r) \Lambda(n+2r): 1 \leq n,r \leq N ) = C_3 + o(1)
\end{equation}
where
\begin{align*}
C_3 &:= \prod_p \frac{\P( n, n+r, n+2r \hbox{ coprime to } p | n, r \in \Z/p\Z )}{
\prod_{j=0}^2 \P( n+jr \hbox{ coprime to } p | n,r \in \Z/p\Z)}\\
&= \frac{3}{2} \prod_{p \geq 5} (1 + \frac{1}{(p-1)^3}) \\
&= 1.534\ldots
\end{align*}
More generally, the Hardy-Littlewood prime tuple conjecture implies that
\begin{equation}\label{hl-ap}
\E( \Lambda(n) \ldots \Lambda(n+(k-1)r): 1 \leq n,r \leq N ) = C_k + o_k(1)
\end{equation}
for all $k \geq 0$ (with the error term $o_k(1)$ depending on $k$), where $C_k$ is the constant
\begin{align*}
C_k &= \prod_p \frac{\P( n, \ldots, n+(k-1)r \hbox{ coprime to } p | n, r \in \Z/p\Z )}{
\prod_{j=0}^{k-1} \P( n+jr \hbox{ coprime to } p | n,r \in \Z/p\Z)}
\end{align*}
which is explicitly computable for each $k$.  The case $k=0$ is trivial, the cases $k=1,2$ follow from the prime number theorem, and
the case $k=3$ is just \eqref{vdc}.  More recently, we have the following results:

\begin{theorem}\label{gt-thm}\cite{gt-primes}, \cite{gt-qm}  The conjecture \eqref{hl-ap} is also true for $k=4$ (so there are infinitely many prime arithmetic progressions of length 4).  Furthermore, for all $k \geq 0$ we have
\begin{equation}\label{lambda-chunk}
\E( \Lambda(n) \ldots \Lambda(n+(k-1)r): 1 \leq n,r \leq N ) > c_k + o_k(1)
\end{equation}
for some explicit constant $c_k > 0$ (which is unfortunately much smaller than $C_k$).  This weaker statement still suffices to establish infinitely many prime arithmetic progressions of length $k$.
\end{theorem}

All of these results have the flavor of ``establish bounds or asymptotics for multilinear averages of $\Lambda$''.  However, some are significantly
harder than others, depending on the exact structure of the multilinear average involved.  As mentioned earlier, the situation has some parallels with the linear, bilinear,
and trilinear Hilbert transform in harmonic analysis; while these expressions are formally very similar in structure, the analytical treatment
of each one in the sequence has proven to be significantly harder than the previous one, for instance no $L^p$ estimates for the trilinear Hilbert transform are currently known.  A certain subclass of these multilinear averages (the ``rank one'' averages involving three or more copies
of $\Lambda$) can be treated by Fourier methods; this includes Vinogradov's theorem and van der Corput's theorem, and see also \cite{balog1} for further discussion.  However, it is by now well established that these techniques cannot directly extend to handle other multilinear averages.  The $k=4$ result in Theorem \ref{gt-thm}
requires a ``quadratic'' generalization of Fourier analysis, pioneered by Gowers \cite{gowers-4}, but still in a very early stage of development.
The higher cases $k \geq 5$ could in principle be treated by polynomial Fourier analysis, of the type developed in \cite{gowers}; this would likely
establish \eqref{hl-ap} for all $k$, this project is currently a work in progress with the author and Ben Green.  Instead, we use an alternate
argument based on ergodic theory which is technically simpler but only gives the weaker result \eqref{lambda-chunk}.

There are two main strategies to obtain progressions:

\begin{itemize}

\item (Uniformity strategy) Attempt to approximate $\Lambda$ by some averaged version $\E(\Lambda|\B)$ of itself, in such a manner that $\Lambda - \E(\Lambda|\B)$ is uniform of the correct order (linearly uniform for $k=3$, quadratically uniform for $k=4$).  This requires one to estimate exponential sums such as $\sum_n \Lambda(n) e(n\theta)$ or $\sum_n \Lambda(n) e(P(n))$ where $P$ is a polynomial or ``generalized polynomial).

\item (Szemer\'edi strategy) Attempt to leverage Szemer\'edi's theorem (or in the case of progressions of length three, Roth's theorem) in order to obtain arithmetic progressions regardless of whether $\Lambda$ is uniform or not.

\end{itemize}

In the case of progressions of length three, the uniformity strategy (more commonly known in this context as the \emph{Hardy-Littlewood circle method})
was developed far earlier than the Szemer\'edi strategy.  It gives sharper results (in particular, it yields the asymptotic \eqref{hl-ap}), but is technically more difficult to implement.  We now briefly discuss each of these strategies in turn.

\section{The uniformity strategy}\label{unif-sec}

We begin by discussing the uniformity strategy.  We shall eschew the traditional framework of the Hardy-Littlewood circle method (which is
only effective for the $k=3$ case) and present this strategy in a language which more easily lends itself to generalization to higher $k$.

The circle method relies on Fourier analysis on the integers $\Z$ (so that the dual group is the unit circle $S^1$, hence the terminology
``circle method'').  For us it will be slightly more convenient to work in the cyclic group $\Z/N\Z$, which is self-dual.  To simplify
the exposition we shall pretend that $\Lambda$ is actually a function on $\Z/N\Z$ rather than $\Z$.  In practice one would have to
justify this by a truncation trick, for instance cutting off $\Lambda$ to $\{1,\ldots,N/3\}$ (possibly using a smooth cutoff function)
and then transferring this to $\Z/N\Z$; this type of ``transference'' is quite standard and introduces no substantial difficulties,
and so we shall gloss over this entire issue.  

Using the above ``cheat'', we can morally rewrite \eqref{hl-ap} as
$$\E( \Lambda(n) \ldots \Lambda(n+(k-1)r): n,r \in \Z/N\Z ) = C_k + o_k(1).$$
Let us first discuss the $k=3$ case (i.e. \eqref{vdc}), which with our new cheat becomes
$$ \E( \Lambda(n) \Lambda(n+r) \Lambda(n+2r): n, r \in \Z/N\Z ) = C_3 + o(1).$$.  
The strategy is to use some variant\footnote{Strictly speaking, one has to replace this Proposition by a weighted variant to cope with the fact that $\Lambda$ is not a bounded function.  This can be done by using a suitable weight function $\Lambda_R$ which is adapted to ``almost primes'', and which among other things obeys a good Fourier restriction theorem which allows one to transfer Proposition \ref{gbgb} to the weighted setting.  See \cite{green}, \cite{green-tao} for further discussion of this issue.} of Proposition \ref{gbgb}.  More specifically, we would seek to approximate $\Lambda$ by an averaged version $\E(\Lambda|\B)$
 such that we have a uniformity estimate\begin{equation}\label{llb}
\| (\Lambda - \E(\Lambda|\B))^\wedge \|_\infty = o(1),
\end{equation}
which (by a suitable variant of Proposition \ref{gbgb}) should imply
\begin{align*} &\E( \Lambda(n) \Lambda(n+r) \Lambda(n+2r): n,r \in \Z/N\Z ) \\
&\quad
= \E( \E(\Lambda|\B)(n) \E(\Lambda|\B)(n+r) \E(\Lambda|\B)(n+2r): n,r \in \Z/N\Z ) + o(1)
\end{align*}
and then one only has to prove \eqref{vdc} for the averaged function $\E(\Lambda|\B)$:
\begin{equation}\label{lam-av}
 \E( \E(\Lambda|\B)(n) \E(\Lambda|\B)(n+r) \E(\Lambda|\B)(n+2r): n,r \in \Z/N\Z ) = C_3 + o(1).
 \end{equation}
 
The first issue is to decide what function $\E(\Lambda|\B)$ to use as the approximant to $\Lambda$.  In order to establish \eqref{lam-av}
we would like $\E(\Lambda|\B)$ to have low ``complexity'' - in particular, it should be far more regular than $\Lambda$ itself - but not so
simple that the approximation to $\Lambda$ is poor in the sense that \eqref{llb} fails.

Let us understand what \eqref{llb} means.  We can rewrite it as
$$ \widehat{\E(\Lambda|\B)}(\xi) = \hat \Lambda(\xi) + o(1) \hbox{ for all } \xi \in \Z/N\Z$$
or in other words
\begin{equation}\label{enxi}
 \E( \E( \Lambda|\B )(n) e_N( -n \xi) | n \in \Z/N\Z ) = \E( \Lambda(n) e_N(-n\xi) | n \in \Z/N\Z ) + o(1)
\hbox{ for all } \xi \in \Z/N\Z.
\end{equation}
This gives us some clues as to what kind of approximation $\E(\Lambda|\B)$ we should choose.  For instance, setting $\xi = 0$
and using the prime number theorem $\E( \Lambda(n) | n \in \Z/N\Z) = 1 + o(1)$, we see that we need $\E(\Lambda|\B)$ to obey
the condition
$$ \E( \E( \Lambda|\B )(n) | n \in \Z/N\Z ) = 1 + o(1).$$
This suggests using the constant function $1$ (or perhaps $\E(\Lambda) = 1+o(1)$) as the approximating function $\E(\Lambda|\B)$; this corresponds 
to interpreting $\B$ as the trivial $\sigma$-algebra $\B_1 := \{ \emptyset, \Z/N\Z\}$.  For this approximation, the left-hand side of
\eqref{lam-av} is very easy to compute, indeed it is just $1 + o(1)$.  Unfortunately, while \eqref{enxi} is true for this approximation
when $\xi = 0$, it is not true for some other values of $\xi$.  Take for instance $\xi = \lfloor N/2 \rfloor$.  Then $e_N(-n\xi)$ is essentially
$+1$ when $n$ is even and $-1$ when $n$ is odd, and so if $\E(\Lambda|\B_1)$ were constant then the left-hand side of \eqref{enxi} would vanish.
On the other hand, the right-hand side of \eqref{enxi} is large and negative, because $\Lambda$ is overwhelmingly supported on the odd numbers
rather than the even numbers.  Thus we must modify the approximant $\E(\Lambda|\B_1)$ to reflect this ``bias'' that $\Lambda$ has towards being
odd.  The easiest way to fix this is to refine the $\sigma$-algebra $\B_1$ to include the odd and even numbers.  In other words, if we now let
$\B_2$ be the $\sigma$-algebra generated by $\B_1$ and the residue classes mod 2 (i.e. the odd and even numbers), then we can use
$\E( \Lambda | \B_2 )$ as our approximant.  By the prime number theorem (and the fact that almost all primes are odd), we know that
this function is $2+o(1)$ on the odd numbers and $o(1)$ on the even numbers.  One can now also check that \eqref{enxi} is now true when $\xi$ is close to zero or close to $N/2$.  Furthermore, the left-hand side of \eqref{lam-av} is quite easy to compute, it is
$$
\frac{\P( n, n+r, n+2r \hbox{ coprime to } 2 | n, r \in \Z/2\Z )}{
\prod_{j=0}^2 \P( n+jr \hbox{ coprime to } 2 | n,r \in \Z/2\Z)} + o(1)
= 2 + o(1).$$
Unfortunately, there are still some further Fourier-analytic biases in $\Lambda$ which are not detected by the approximation
$\E( \Lambda | \B_2 )$, for instance the fact that $\Lambda$ is concentrated in the residue classes $1 \mod 3$ and $2 \mod 3$ and nearly
vanishes on the residue class $0 \mod 3$ will cause the Fourier coefficients of $\Lambda$ to be rather large for $\xi$ near $N/3$ and $2N/3$, whereas
$\E( \Lambda | \B_2 )$ is uniformly distributed among all three residue classes $\mod 3$ and thus has a negligible Fourier coefficient
at those frequencies.  One can address this failure of \eqref{enxi} by refining the approximation $\E( \Lambda | \B_2 )$ further to $\E( \Lambda | \B_3 )$, where $\B_3$ is the $\sigma$-algebra formed by adjoining the residue classes modulo 3 to $\B_2$ (or in other words, $\B_3$ is the $\sigma$-algebra
generated by the residue classes modulo 6).  Then one can show that \eqref{enxi} now holds for all $\xi$ near multiples of $N/6$.  Furthermore,
one has $\E( \Lambda | \B_3 )(n) = 3+o(1)$ when $n$ is coprime to 6 and $\E( \Lambda | \B_3 ) = o(1)$ otherwise; this follows from the
prime number theorem combined with \emph{Dirichlet's theorem}, which asserts that $\Lambda$ is uniformly distributed among those
residue classes modulo $m$ which are coprime to $m$, as long as $N$ is sufficiently large compared to $m$ (here we take $m=6$).
Because of this, one can compute (using the Chinese remainder theorem) that the left-hand side of \eqref{lam-av} is now 
$$
\prod_{p=2,3} \frac{\P( n, n+r, n+2r \hbox{ coprime to } p | n, r \in \Z/2\Z )}{
\prod_{j=0}^2 \P( n+jr \hbox{ coprime to } p | n,r \in \Z/2\Z)}
+ o(1) = \frac{3}{2}  + o(1).$$

One can of course continue in this fashion.  Let $w = w(N)$ be a slowly growing function of $N$, e.g. $w = \log \log N$, and let $W$ be the
product of all the primes less than $w$.  We let $\B_w$ be the $\sigma$-algebra formed by the residue classes modulo $W$, then we use
$\E(\Lambda|\B_w)$ as our approximant.  From Dirichlet's theorem, one can show (if $w$ is sufficiently slowly growing in $N$) that
$\E(\Lambda|\B_w)(n) = \frac{W}{\phi(W)} + o(1)$ if $n$ is coprime to $W$, and $\E(\Lambda|\B_w)(n) = o(1)$ otherwise; here $\phi(W)$
is the Euler totient function of $W$, i.e. the number of integers in $\{1,\ldots,W\}$ which are coprime to $W$.  From the Chinese
remainder theorem, the left-hand side of \eqref{lam-av} can be computed as
$$
\prod_{p \leq w} \frac{\P( n, n+r, n+2r \hbox{ coprime to } p | n, r \in \Z/2\Z )}{
\prod_{j=0}^2 \P( n+jr \hbox{ coprime to } p | n,r \in \Z/2\Z)}
+ o(1) = C_3 + o(1)$$
since the product is convergent and $w$ tends (slowly) to infinity.  Thus it only remains to demonstrate \eqref{enxi}.  This would be
easy if $w$ was extremely large (e.g. if $w = \sqrt{N}$, then the sieve of Eratosthenes essentially ensures that
$\Lambda = \E( \Lambda|\B_w )$, but unfortunately the error terms blow up long before $w$ reaches this level.  Nevertheless,
this ``$W$-trick'' of removing all the structure from $\Lambda$ associated to those primes less than $w$ does make the task
of \eqref{enxi} much easier.  Essentially, it means that \eqref{enxi} is automatically true whenever $\xi$ is a ``major arc frequency'',
which roughly means that $\xi \approx aN/q$ for some integers $a,q$ with $q \leq w$.  It thus remains to prove \eqref{enxi} when $\xi$
is a ``minor arc'' frequency, which roughly means that $q\xi$ is not close to zero modulo $N$ for any $q \leq w$.  In such a case,
the left-hand side of \eqref{enxi} is very small (by construction of $\E(\Lambda|\B_w)$, and one is reduced to establishing enough cancellation
in the sum $\sum_{n < N} \Lambda(n) e_N(-n\xi)$ to ensure that it is $o(N)$.  (Note that the trivial bound coming from using
absolute values and the prime number theorem is $O(N)$).

To do this, one must finally use some deeper structure of the function $\Lambda(n)$, beyond the prime number theorem
and Dirichlet's theorem.  This was first done by Vinogradov, with later simplifications by Vaughan and other authors; we present a vastly oversimplified
sketch of the main idea here.  The starting point is the identity \eqref{logn}.  Solving for $n$ we obtain
the formula
$$ \Lambda(n) = \sum_{c,d: cd = n} \log c \mu(d),$$
where $\mu(d)$ is the \emph{M\"obius function}, defined as $(-1)^m$ if $d$ is the product of $m$ distinct primes, and equal to 0 otherwise.
Thus we can write
$$ \sum_{n < N} \Lambda(n) e_N(-n\xi) = \sum_{c,d: cd < N} \log c \mu(d) e_N(-cd\xi).$$
The idea is now to view this as a bilinear form acting on the functions $\log$ and $\mu$, given by the matrix coefficients $e_N(-cd\xi)$.  The
hypothesis that $\xi$ is not ``minor arc'' leads to some almost orthogonality in this matrix (which can be made explicit by the $TT^*$ method),
which after some care can eventually lead to the $o(1)$ gain.  (This is an oversimplification because the portions of this expression when
$c$ or $d$ is small require some additional attention, including a quantitative version of Dirichlet's theorem known as the Siegel-Walfisz theorem; we will not discuss these rather lengthy issues here).  This can eventually be used to establish Van der Corput's theorem \eqref{vdc}.

It turns out that the same ideas can also be pushed (with several additional difficulties) to give the $k=4$ case of \eqref{hl-ap}; it is not 
yet known whether the arguments can be pushed to general $k$.  By using a result similar to Lemma \ref{gvn-thm}, as a substitute for Proposition \ref{gbgb}, it suffices to find an approximation $\E(\Lambda|\B)$ for $\Lambda$ such that
\begin{equation}\label{lab-4}
 \E( \E(\Lambda|\B)(n) \E(\Lambda|\B)(n+r) \E(\Lambda|\B)(n+2r) \E(\Lambda|\B)(n+3r): n,r \in \Z/N\Z ) = C_4 + o(1)
 \end{equation}
and
\begin{equation}\label{u3}
 \| \Lambda - \E(\Lambda|\B) \|_{U^3} = o(1).
 \end{equation}
As in the $k=3$ case, we again invoke the ``$W$-trick'' and set $\B = \B_w$ where $w$ is again a slowly growing function of $N$.  
When one does so, \eqref{lab-4} is easy to establish, but \eqref{u3} is still quite difficult.  Expanding out the $U^3$ norm directly
gives rise to expressions which are about as complicated to estimate as the original expression in \eqref{hl-ap}.  However, one can proceed instead
by using the inverse theory used in Gowers' proof of Szemer\'edi's theorem for progressions of length 4.  The idea is to assume that
$\| \Lambda - \E(\Lambda|\B) \|_{U^3}$ is large, say larger than some $\delta > 0$, and arrive at a contradiction.  One can repeat the analysis
in Gowers' arguments (though one has to introduce weights to deal with the fact that $\Lambda$ is not bounded) to eventually conclude that
$$ \E( (\Lambda(n) - \E(\Lambda|\B)(n)) e_N( P(n) ) | n \in \Z/N\Z ) \geq c(\delta) > 0$$
for some ``generalized quadratic phase function'' $P(n)$; we shall gloss over exactly what ``generalized quadratic phase function'' means
here but one should think of $P$ as being like a quadratic polynomial.  Thus to conclude the proof, one needs to extend the linear uniformity
estimate \eqref{enxi} to the claim that
$$  \E( \E( \Lambda|\B )(n) e_N( P(n)) | n \in \Z/N\Z ) = \E( \Lambda(x) e_N(P(n)) | n \in \Z/N\Z ) + o(1)$$
for all generalized quadratic phase functions $P$.  It turns out that once again one can divide into the case when $P$ is ``major arc'' - 
all the non-constant coefficients of $P$ are essentially rational multiples of $N$ with small denominator, and when $P$ is ``minor arc'' -
when at least one of the coefficients behaves ``irrationally''.  The major arc case is again easy, while the minor arc case turns out
to be again amenable to the methods of Vinogradov and Vaughan.  Here the point is to establish some orthogonality in 
the matrix coefficients $e_N(P(cd))$.  See \cite{gt-qm} for further details.

\section{The Szemer\'edi strategy}

In principle, the uniformity strategy discussed above should in fact prove \eqref{hl-ap} for all $k$.  However, at present we are
restricted to $k \leq 4$ because the inverse theorem that passes from large $U^{k-1}$ norm to correlation with a generalized polynomial
phase function of order $k-2$ has only been rigorously proven for $k \leq 4$.  (The analysis in \cite{gowers} strongly suggests that this inverse theorem should in fact extend to higher $k$; this is a current work in progress with the author and Ben Green).  In particular,
while it is conjectured that we in fact have
\begin{equation}\label{bbw}
 \| \Lambda - \E(\Lambda|\B_w) \|_{U^{k-1}} = o_k(1)
\end{equation}
for all $k$ (which would certainly imply \eqref{hl-ap}), this estimate has not yet been rigorously established.

Nevertheless, one can still achieve the weaker statement \eqref{lambda-chunk} by using ergodic theory arguments to locate 
another $\sigma$-algebra $\B$ (which could be somewhat finer than $\B_w$) for which the analogue of \eqref{bbw} holds.  To finish the proof of \eqref{lambda-chunk}, it then remains to show that
\begin{equation}\label{lab}
\E( \E(\Lambda|\B)(n) \ldots \E(\Lambda|\B)(n+(k-1)r): 1 \leq n,r \leq N ) > c_k + o_k(1).
\end{equation}
Unfortunately, the structure of the algebra $\B$ is much less well understood than $\B_w$, and as such the function $\E(\Lambda|\B)$ is also
not very well understood.  However, being a conditional expectation of $\Lambda$, it is still non-negative, has the same mean (i.e. $1 + o(1)$)
as $\Lambda$.  Crucially, one can also establish that $\E(\Lambda|\B)$ is also \emph{bounded} by $O(1)$.  By the third version of Szemer\'edi's theorem,
these three facts imply \eqref{lab}.

A prototype of this argument is the proof of Theorem \ref{gthm} in \cite{green}, which used Fourier analytic methods (but with ergodic ideas lurking under the surface), and as such was limited to the $k=3$ case.  This argument was then simplified and extended in \cite{green-tao}; simultaneously, in \cite{gt-primes} the Fourier-analytic components were replaced with ergodic theory arguments which could then extend to general $k$.
Here we shall begin by discussing the general ergodic theory argument, and return to briefly discuss the earlier Fourier-analytic arguments at
the end of this section.

One important technical problem
that needs addressing is that the function $\Lambda$ is not bounded, which means that much of the analysis in previous sections, strictly speaking,
does not apply.  This is essentially equivalent to the fact that the primes have asymptotic density zero.
However, one can resolve this problem by bounding $\Lambda$ not by a bounded multiple of the constant function 1, which is not possible, but instead by a bounded multiple of another function $\nu$ which resembles $\Lambda$ but is much easier to work with\footnote{As before we are ignoring some details concerning how one embeds $\Lambda$ inside $\Z/N\Z$; also, it turns out to be convenient to ``factor out'' the initial $\sigma$-algebra $\B_w$ by passing to a single atom, such as the residue class $1 \mod W$; we ignore these minor technical issues here.}.  This corresponds to 
viewing the primes not as a
(sparse) subset of the integers, but rather as a subset of the set of \emph{almost primes}, which is much more tractable than the primes to study,
and with the property that the primes have positive \emph{relative} density inside the primes.  One byproduct of this approach is that,
because it uses very little about the primes other than this positive relative density, it in fact implies a stronger result, namely that
\emph{all} subsets of the primes with positive relative density must necessarily contain arbitrarily long arithmetic progressions.

Informally, the idea is as follows.  Let $P$ be the set of prime numbers between $N/2$ and $N$. The sieve of Eratosthenes shows that $P$ consists precisely of those integers in $\{N/2, \ldots, N\}$ which are coprime to all primes less than $\sqrt{N}$.  Motivated by this, let us define the partially sifted set $P_R$ to be those integers in $\{N/2, \ldots, N\}$ which are coprime to all primes less than $R$, where $1 \leq R \leq \sqrt{N}$
is a parameter.  Thus as $R$ increases to $\sqrt{N}$, $P_R$ decreases until it becomes $P$.  The first few sets $P_R$ are easy to understand,
for instance $P_2$ is simply the odd numbers from $N/2$ to $N$.  In particular, any statistic involving $P_R$ (e.g. counting how many arithmetic
progressions of length $k$ are contained in $P_R$) is quite easy to compute to high accuracy when $R$ is small.  However, the task becomes
increasingly difficult when $R$ gets large.  The vast and well-developed topic of \emph{sieve theory} - a key component of analytic number theory -
is devoted to questions like this; while this theory is too complex to be surveyed here, let us oversimplify one of the basic results in
that field, namely the \emph{fundamental lemma of sieve theory}  In our notation, this lemma roughly speaking asserts that
that one can compute the statistics of $P_R$ as long as $R$ is a sufficiently small power of $N$.  For instance, one can accurately count
the number of arithmetic progressions in $P_R$ of length $k$ if $R$ is less than $N^{1/2k}$.

An informal probabilistic argument suggests that
$$ |P_R| \sim \frac{N}{2} \prod_{p < R} (1 - \frac{1}{p}) \sim \frac{N}{\log R}$$
where we use $X \sim Y$ to denote equivalence up to constants (i.e. $C^{-1} Y \leq X \leq CY$).  A famous theorem of Merten in fact gives the
more precise asymptotic
$$ |P_R| = (1 + o(1)) \frac{N}{2} \frac{e^{-\gamma}}{\log R}$$
as long as $R$ is much less than $\sqrt{N}$ but goes to infinity as $N \to \infty$.  Here $\gamma = 0.577\ldots$ is Euler's constant.
Comparing this with the prime number theorem
$$ |P| = (1 + o(1)) \frac{N}{2} \frac{1}{\log N}$$
we see that $P$ will have a relative density $|P|/|P_R|$ bounded away from zero as long as we set $R$ to equal a small power of $N$, say
$R = N^{\eps}$ for some fixed $\eps$ (this $\eps$ will eventually depend on $k$; in \cite{gt-primes} it is $\eps = \frac{1}{k 2^{k+4}}$).

A natural choice for the weight function $\nu$ would then be $\log R 1_{P_R}$; this function would thus be normalized to
essentially have mean 1, and $\Lambda$ would be dominated by a bounded multiple of $\nu$.  For technical reasons, however, the function
$1_{P_R}$ is a bit too ``rough'' to serve as a good weight function, and it is better to use a slightly ``smoother'' variant of this
function, namely the truncated divisor sums studied by Goldston and Yildirim \cite{goldston-yildirim-old1}, \cite{goldston-yildirim-old2},
\cite{goldston-yildirim}.  These are formed by replacing the von Mangoldt function
$$ \Lambda(n) = \sum_{d|n} \mu(d) \log \frac{n}{d}$$
with the variant
$$ \Lambda_R(n) := \sum_{d|n} \mu(d) (\log \frac{R}{d})_+$$
where $x_+ := \max(x,0)$ is the positive part of $x$.  One can easily verify that $\Lambda_R$ is equal to $\log R$ on the set of
$P_R$, and can thus be thought of as the function $\log R 1_{P_R}$ with an additional ``tail''.  The advantage of working with $\Lambda_R$
instead of $P_R$ is that $\Lambda_R$ is easily expressed as a linear combination of the functions $1_{d|n}$, i.e. the characteristic functions of
the residue class $d \Z$.  Moreover, the coefficients $\mu(d) (\log \frac{R}{d})_+$ for this linear combination are supported on the small
values of $d$, which are easier to control; this is roughly analogous in harmonic analysis to a function having Fourier transform supported on
the ``low frequencies'', which explains why such functions in number theory are sometimes referred to as being ``smooth''.  In particular, the work
of Goldston and Yildirim showed that (providing $R$ was a sufficiently small power of $N$) it was possible to accurately estimate such expressions as
$$ \E( \Lambda_R(n) \Lambda_R(n+r) \ldots \Lambda_R(n+(k-1)r) | n, r \in \Z/N\Z ).$$
We cannot directly use $\Lambda_R$ to dominate $\Lambda$, as it turns out to oscillate in sign; however this is easily fixed by using instead
the function $\nu(n) := \frac{1}{\log R} \Lambda_R^2(n)$.  Actually, this is an oversimplification; in practice we need to localize $n$ to
an arithmetic progression of spacing $W$ and length equal to $cN$ for a small multiple of $N$.  After these adjustments, Goldston and Yildirim
essentially showed that $\nu$ was ``pseudorandom'' - that almost all the correlations of $\nu$ were very close to 1 (a formal definition of this rather technical statement is in \cite{gt-primes}).  Another way of saying this is that $\nu$ lies very close to 1 in certain ``weak'' norms (such as
the Gowers uniformity norms).  With this pseudorandomness property, it turns out that the weight $\nu$ behaves very similarly to 1, thus for instance
the generalized von Neumann theorem, Lemma \ref{gvn-thm}, can be extended to the case where $f$ is bounded by the pseudorandom function $\nu$
rather than the constant function $1$ (although one has to accept some additional $o(1)$ errors when doing so).  See \cite{gt-primes} for details;
the ideas here were initially motivated by similar arguments in the setting of hypergraphs by Gowers \cite{gowers-reg}.

We can now describe the proof of \eqref{lambda-chunk} for general $k$.  For sake of concreteness we shall restrict ourselves to the case $k=4$,
although the argument extends without difficulty to higher $k$.  We shall use the machinery developed in the energy increment proof of Szemer\'edi's theorem in the $k=4$ case.

As discussed earlier, the objective is to locate a $\sigma$-algebra $\B$ such that 
\begin{equation}\label{u-small}
\| \Lambda - \E(\Lambda|\B) \|_{U^{k-1}} \hbox{ is small}
\end{equation}
(where we shall be a bit vague as to what ``small'' means), and such that $\E( \Lambda | \B )$ is bounded.  The choice $\B = \B_w$, where $w$
as before is a slowly growing function of $N$, will obey the second property (this is basically Dirichlet's theorem), but it is unknown
as to whether it obeys the first property.  Nevertheless, we can proceed by a stopping time argument, somewhat similar to the Calder\'on-Zygmund
stopping time arguments used in harmonic analysis, or the stopping time argument used in the proof of the Szemer\'edi regularity lemma.  The key
point is that if \eqref{u-small} fails for some algebra $\B$, then by setting $g$ to be the dual function of $\Lambda - \E(\Lambda|\B)$,
$$ g := {\mathcal D}(\Lambda - \E(\Lambda|\B)),$$
then $g$ will have a non-trivial correlation with $\Lambda - \E( \Lambda | \B )$:
$$ |\langle g, \Lambda - \E( \Lambda | \B ) \rangle| \hbox{ is large}.$$
Viewing this geometrically in the Hilbert space $L^2(\Z_N)$, this means that $\Lambda$ (now thought of as a vector) contains a non-trivial 
component which is orthogonal to the subspace $L^2(\B)$ which the conditional expectation operator $\E(|\B)$ projects to, and which is also
somewhat parallel to $g$.  Thus if one defines $\B'$ to be the algebra generated by $\B$ and (suitable level sets of) $g$, we expect
$L^2(\B')$ to capture both $L^2(\B)$ and $g$ (or a vector very close to $g$).  Putting this together, we expect $\Lambda$ to be closer to
the subspace $L^2(\B')$ than to the smaller subspace $L^2(\B)$; indeed, some applications of Cauchy-Schwarz and Pythagoras's theorem can be
used to give an energy increment estimate of the form
\begin{equation}\label{energy-increment}
\| \E( \Lambda | \B' ) \|_{L^2}^2 \geq \| \E( \Lambda | \B ) \|_{L^2}^2 + c
\end{equation}
for some $c > 0$ (which depends of course on the definitions of ``small'' and ``large'').  

To summarize, whenever \eqref{u-small} fails, we can exploit this failure to enlarge the underlying $\sigma$-algebra $\B$ in such a way that
it collects more of the ``energy'' of $\Lambda$.  We can now replace $\B$ by $\B'$ and iterate this procedure until \eqref{u-small} is finally
attained.  At first glance it seems that this algorithm could continue for quite a long time, since $\Lambda$ has a large $L^2$ norm.
Fortunately, though, it turns out that $\E( \Lambda | \B )$ remains uniformly bounded throughout this algorithm.  This is because $\Lambda$
is bounded by $\nu$, and thus $\E( \Lambda | \B )$ is bounded by $\E( \nu | \B )$.  The latter function turns out to be bounded because $\nu$ is pseudorandom (and thus very uniform), whereas $\B$ was essentially generated by dual functions (and thus highly non-uniform).  Indeed,
it turns out that even if one runs this algorithm for a large number of iterations, the bounds on $\E( \nu | \B )$ only worsen by at most
$o(1)$.  This crucial fact is one of the more delicate computations in \cite{gt-primes}, but it ultimately follows from
the pseudorandomness information on $\nu$ and an application of the Gowers-Cauchy-Schwarz inequality \eqref{gcz}.
This boundedness of $\E( \Lambda | \B )$ is required for two reasons: firstly, in order that Szemer\'edi's theorem (in its third
formulation) can be applied to this function, and secondly it is used (in conjunction with \eqref{energy-increment}) to show that
the algorithm to find $\B$ halts after only a bounded number of iterations.

We now briefly remark on the earlier $k=3$ versions of the above argument, referring the reader to \cite{green}, \cite{green-tao} for further details.  In that case, the notion of pseudorandomness of the dominating measure $\nu$ was replaced by that of \emph{linear pseudorandomness} or \emph{Fourier pseudorandomness}, which basically asserts that all the Fourier coefficients of $\nu-1$ were small.  By Tomas-Stein restriction type arguments, this implies a certain Fourier restriction theorem for $\nu$, which can be used to develop weighted analogues of Proposition \ref{gbgb} adapted to $\nu$.
One then runs the same argument as before, but this time the $\sigma$-algebra $\B$ is more explicit: it is the algebra generated by the Bohr sets corresponding to those frequencies where the Fourier transform of $\Lambda$ is large.  (Of course, the Hardy-Littlewood method already
provides information as to where this Fourier transform is large; however the advantage of this argument is that it still works if $\Lambda$ is
replaced by any other function supported on a dense subset of the primes, whereas the Hardy-Littlewood method relies on
the arithmetic structure on $\Lambda$ and does not extend in this manner).  Again, the pseudorandomness of $\nu$ will ensure that $\E( \nu | \B )$,
and hence $\E( \Lambda | \B )$, is bounded, and one can then apply (the third version of) Roth's theorem to deduce Theorem \ref{gthm}.
(Some further variations of this theme are pursued in \cite{green-tao}).

\providecommand{\bysame}{\leavevmode\hbox to3em{\hrulefill}\thinspace}
\providecommand{\MR}{\relax\ifhmode\unskip\space\fi MR }
\providecommand{\MRhref}[2]{%
  \href{http://www.ams.org/mathscinet-getitem?mr=#1}{#2}
}
\providecommand{\href}[2]{#2}


\begin{thebibliography}{10}

\bibitem{assani}
I. Assani, \emph{Pointwise convergence of ergodic averages along cubes}, preprint.

\bibitem{balog1} A. Balog, \emph{Linear equations in primes,} Mathematika \textbf{39} (1992) 367--378.

\bibitem{bergelson-leibman} V. Bergelson and A. Leibman, \emph{Polynomial extensions of van der Waerden's and Szemer\'edi's theorems,} J. Amer. Math. Soc. \textbf{9} (1996), 725--753.

\bibitem{bourg2} J. Bourgain \emph{A Szemer\'edi-type theorem for sets of positive density in $\mathbb{R}^k$,} Israel J. Math \textbf{54} (1986), no. 3, 307--316.

\bibitem{bourgain-triples}
\bysame, \emph{On triples in arithmetic progression}, GAFA \textbf{9} (1999), 968--984.

\bibitem{chang}
F.R.K. Chang, R. Graham, \emph{Quasi-random subsets of $\Z_n$}, J. Comb. Th. A. \textbf{61} (1992), 64--86.

\bibitem{chowla}
S. Chowla, \emph{There exists an infinity of 3---combinations of primes in A. P.},
 Proc. Lahore Philos. Soc. \textbf{6}, (1944). no. 2, 15--16.

\bibitem{ccw}
M. Christ, A. Carbery, J. Wright, \emph{Multidimensional van der Corput and Sublevel set estimates}, J. Amer. Math. Soc. \textbf{12} (1999), 981-1015.

\bibitem{erdos}
P. Erd\H{o}s, P. Tur\'an, \emph{On some sequences of integers}, J. London Math. Soc. \textbf{11} (1936), 261--264.

\bibitem{furst}
H. Furstenberg, \emph{Ergodic behavior of diagonal measures and a theorem of Szemer\'edi on arithmetic progressions}, J. Analyse Math. \textbf{31} (1977), 204--256.

\bibitem{FKO}
H. Furstenberg, Y. Katznelson and D. Ornstein, \emph{The ergodic-theoretical proof of Szemer\'edi's theorem,} Bull. Amer. Math. Soc. \textbf{7} (1982), 527--552.

\bibitem{furst-weiss}
H. Furstenberg, B. Weiss, \emph{A mean ergodic theorem for $1/N \sum_{n=1}^N f(T^n x) g(T^{n^2} x)$}, Convergence in ergodic theory and probability (Columbus OH 1993), 193--227, Ohio State Univ. Math. Res. Inst. Publ., 5. de Gruyter, Berlin, 1996.

\bibitem{goldston-yildirim-old1} D. Goldston and C.Y. Yildirim \emph{Higher correlations of divisor sums related to primes, I: Triple correlations,} Integers \textbf{3} (2003) A5, 66pp.

\bibitem{goldston-yildirim-old2} \bysame, \emph{Higher correlations of divisor sums related to primes, III: $k$-correlations,} preprint (available at AIM preprints)

\bibitem{goldston-yildirim} \bysame, \emph{Small gaps between primes, I,} preprint.

\bibitem{gowers-4}
T. Gowers, \emph{A new proof of Szemer\'edi's theorem for arithmetic 
progressions of length four}, GAFA \textbf{8} (1998), 529--551.
 
\bibitem{gowers}
\bysame, \emph{A new proof of Szemer\'edi's theorem}, GAFA \textbf{11} (2001), 465-588.

\bibitem{gowers-reg}
\bysame, \emph{Hypergraph regularity and the multidimensional Szemer\'edi theorem,} preprint

\bibitem{graham}
R. Graham, B. Rothschild, J.H. Spencer, Ramsey Theory, John Wiley and Sons, NY (1980).

\bibitem{green}
B.J. Green, \emph{Roth's theorem in the primes,} preprint.

\bibitem{green-reg}
\bysame, \emph{A Szemer\'edi-type regularity lemma in abelian groups,} preprint.

\bibitem{green-survey}
\bysame, \emph{Finite field models in arithmetic combinatorics}, preprint.

\bibitem{gt-primes}
B.J. Green and T. Tao, \emph{The primes contain arbitrarily long arithmetic progressions,} preprint.

\bibitem{green-tao}
B.J. Green and T. Tao, \emph{Restriction theory of Selberg's sieve, with applications,} preprint.

\bibitem{gt-qm}
B.J. Green and T. Tao, \emph{An inverse theorem for the Gowers $U^3$ norm}, preprint.

\bibitem{hardy-littlewood} G.H. Hardy and J.E. Littlewood \emph{Some problems of ``partitio numerorum''; III: On the expression of a number as a sum of primes,} Acta Math. \textbf{44} (1923), 1--70

\bibitem{heath-brown1}
D.R. Heath-Brown, \emph{Three primes and an almost prime in arithmetic progression,} J. London Math. Soc. (2) \textbf{23} (1981), 396--414.

\bibitem{heath-brown2}
\bysame, \emph{Linear relations amongst sums of two squares,} Number theory and algebraic geometry --- to Peter Swinnerton-Dyer on his 75th birthday, CUP (2003). 

\bibitem{host-kra2}
\bysame, \emph{Non-conventional ergodic averages and nilmanifolds,} to appear in Ann. Math. 

\bibitem{laceyt1}M. Lacey, C. Thiele, {\it 
$L^p$ estimates on the bilinear Hilbert transform for $2<p<\infty$.}
Ann. Math. 146 (1997), pp. 693--724.


\bibitem{moran-pritchard-thyssen} A. Moran, P. Pritchard and A. Thyssen, \emph{Twenty-two primes in arithmetic progression,} Math. Comp. \textbf{64} (1995), no. 211, 1337--1339.

\bibitem{ramare} O. Ramar\'e, \emph{On Snirel'man's constant,} Ann. Scu. Norm. Pisa \textbf{21} (1995), 645--706.

\bibitem{ramare-ruzsa} O. Ramar\'e and I.Z. Ruzsa, \emph{Additive properties of dense subsets of sifted sequences,} J. Th. Nombres de Bordeaux \textbf{13} (2001) 559--581.

\bibitem{rodl}
V. R\"odl, J. Skokan, \emph{Regularity lemma for $k$-uniform hypergraphs}, to appear, Random Structures and Algorithms.

\bibitem{rodl2}
V. R\"odl, J. Skokan, \emph{Applications of the regularity lemma for uniform hypergraphs}, preprint.

\bibitem{roth}
K.F. Roth, \emph{On certain sets of integers}, J. London Math. Soc. \textbf{28} (1953), 245-252.

\bibitem{solymosi}
J. Solymosi, \emph{A note on a question of Erd\"os and Graham}, Combinatorics, Probability and Computing \textbf{13} (2004), 263--267.

\bibitem{szemeredi-4}
E. Szemer\'edi, \emph{On sets of integers containing no four elements in arithmetic progression},
Acta Math. Acad. Sci. Hungar. \textbf{20} (1969), 89--104.

\bibitem{szemeredi}
\bysame, \emph{On sets of integers containing no $k$ elements in arithmetic progression},
Acta Arith. \textbf{27} (1975), 299--345.

\bibitem{szem-reg}
\bysame, \emph{Regular partitions of graphs,} in ``Proc. Colloque Inter. CNRS'' (J.-C. Bermond, J.-C. Fournier, M. Las Vergnas, D. Sotteau, eds.) (1978), 399--401.

\bibitem{tao:szemeredi}
T. Tao, \emph{A quantitative ergodic theory proof of Szemer\'edi's theorem}, preprint.

\bibitem{titchmarsh}
E.C. Titchmarsh, \emph{The theory of the Riemann zeta function,} Oxford University Press, 2nd ed, 1986.

\bibitem{van-der-corput}
J.G. van der Corput, \emph{\"Uber Summen von Primzahlen und Primzahlquadraten,} Math. Ann. \textbf{116} (1939), 1--50.

\bibitem{vdw}
B.L. Van der Waerden, \emph{Beweis einer Baudetschen Vermutung}, Nieuw. Arch. Wisk. \textbf{15} (1927), 212--216.

\bibitem{varnavides} P. Varnavides, \emph{On certain sets of positive density,} J. London Math. Soc. \textbf{34} (1959) 358--360.

\bibitem{vinogradov}
I.M. Vinogradov, \emph{Representation of an Odd Number as a Sum of Three Primes}, Comptes rendus (Doklady) de l'Acad\'emie des Sciences de l'U.R.S.S. \textbf{15} (1937a), 169--172. 

\bibitem{ziegler}
T. Ziegler, \emph{Universal characteristic factors and Furstenberg averages}, preprint.

\bibitem{ziegler2}
\bysame, \emph{A non-conventional ergodic theorem for a nilsystem}, preprint.

\end{thebibliography}
     \end{document}